\RequirePackage{fix-cm}
\documentclass[twocolumn]{svjour3}          
\smartqed  
\usepackage{lineno,hyperref}
\usepackage{amssymb}
\usepackage{graphicx}
\usepackage{listings}
\usepackage{color}
\usepackage{float}
\usepackage{subfig}
\usepackage{algorithm}
\usepackage{algorithmic}
\usepackage{amsmath}

\usepackage{enumitem}

\usepackage{booktabs} 

\usepackage{etoolbox}
\newcommand*{\affaddr}[1]{#1} 
\newcommand*{\affmark}[1][*]{\textsuperscript{#1}}
\def\CC{{C\hspace{-.05em}\raisebox{.4ex}{\tiny\bf ++}}}

\newif\ifproofreadnew
\newcommand{\markline}[1]{%
\ifproofreadnew
\textcolor{red}{#1}%
\else
#1%
\fi
}

\newif\ifproofread
\newcommand{\marklineNew}[1]{%
\ifproofread
\textcolor{blue}{#1}%
\else
#1%
\fi
}

\newif\ifproofreadL
\newcommand{\marklineNewL}[1]{%
\ifproofreadL
\textcolor{magenta}{#1}%
\else
#1%
\fi
}

\begin{document}
\proofreadLfalse 
\proofreadfalse 
\proofreadnewfalse 

\title{Topology optimization of transient vibroacoustic problems for broadband filter design using cut elements
}


\author{Cetin B. Dilgen\affmark[*]         \and
		Niels Aage  \thanks{\affaddr{\affmark[*]Corresponding author} }                  
}


\institute{Cetin B. Dilgen \at
Centre for Acoustic-Mechanical Micro Systems (CAMM) \\ Department of Civil and Mechanical Engineering \\ Technical University of Denmark \\ Nils Koppels All\'e, Building 404 \\ DK-2800, Denmark \\
              \email{cedil@mek.dtu.dk}           
           \and
           Niels Aage \at
Centre for Acoustic-Mechanical Micro Systems (CAMM) \\ Civil and Department of Mechanical Engineering \\ Technical University of Denmark \\ Nils Koppels All\'e, Building 404 \\ DK-2800, Denmark \\
		\email{naage@mek.dtu.dk} 
}

\date{Received: date / Accepted: date}

\authorrunning{Cetin B. Dilgen et al.}
\maketitle

\begin{abstract}
The focus of this article is on shape and topology optimization of transient vibroacoustic problems. The main contribution is a transient problem formulation that enables optimization over wide ranges of frequencies with complex signals, which are often of interest in industry. The work employs time domain methods to realize wide band optimization in the frequency domain. To this end, the objective function is defined in frequency domain where the frequency response of the system is obtained through a fast Fourier transform (FFT) algorithm on the transient response of the system. The work utilizes a parametric level set approach to implicitly define the geometry in which the zero level describes the interface between acoustic and structural domains. A cut element method is used to capture the geometry on a fixed background mesh through utilization of a special integration scheme that accurately resolves the interface. This allows for accurate solutions to strongly coupled vibroacoustic systems without having to re-mesh at each design update. The present work relies on efficient gradient based optimizers where the discrete adjoint method is used to calculate the sensitivities of objective and constraint functions. A thorough explanation of the consistent sensitivity calculation is given involving the FFT operation needed to define the objective function in frequency domain. Finally, the developed framework is applied to various vibroacoustic filter designs and the optimization results are verified using commercial finite element software with a steady state time-harmonic formulation.
\keywords{Vibroacoustics \and Cut finite elements \and Immersed boundary methods  \and Transient optimization \and Level set methods \and Shape optimization}
\end{abstract}

\section{Introduction}
Topology optimization \cite{bends2003a} is a numerical method to determine optimal material distributions that minimizes  a given performance criteria under a set of constraints. The method has gained  increasing popularity across many areas in both research and industry since it allows for 
innovative designs through free material distribution. Considerable number of studies on topology optimization paved the way for giga scale structural optimization for static structural mechanics problems \cite{aage2017a}. The method has also been applied to the optimization of fluid systems \cite{kontoleontos2013a,dilgen2018a} as well as several multi-physics applications \cite{sigmund2001a,alexandersen2016a,dilgen2018b}. Considering dynamics, after the introduction of the adjoint method for transient problems \cite{haug1978a,michaleris1994a}
examples of topology optimization under transient loadings can be found in the works of \cite{pedersen2004a,li2004a,turteltaub2005a}. Moreover, optimizing for complex signals containing wide range of frequency content is of great interest in industry, for example, optimizing parts of a hearing aid device which requires a good performance expectancy over a wide frequency window. Commonly, frequency domain modelling and optimization only deals with discrete frequencies which are usually selected from a limited window of the frequency band. Optimization can quickly become excessively expensive when many frequencies are considered in order to widen the zone of influence of optimization on the frequency response of the system. Examples of topology optimization in frequency domain can be found considering eigenvalue problems \cite{seyranian1994a}, optics \cite{jensen2011a}, acoustics \cite{duehring2008a,park2008a,christiansen2016a} and vibroacoustics \cite{yoon2007a}. Using a time-dependent problem formulation presents a promising alternative to address this issue since the optimization can be carried out using complex and compact signals that contains broad ranges of frequencies. This idea can be seen in the works of \cite{nomura2007a,hassan2014a,hassan2015a,Hyun2021} 
where a representative time-domain input pulse is selected to excite a broad frequency range 
in order to carry out transient topology optimization of antennas. 
However, this approach still does not provide the full control on the broad-band response in the frequency domain. This is because the objective function is defined in the time domain and the optimization only indirectly effects the frequency content of the signal that is being optimized.

Although density based topology optimization methods provides the largest degree of design freedom, the ersatz material model presents several issues for multi-phsysics problems that are strongly coupled through the interface. The main reasons for this can be listed as the lack of physical interpretation of the intermediate densities and the stair-case (pixelized) boundary description. Generally, strongly coupled problems, i.e vibroacoustics, require accurate modelling of the interface in order to correctly capture and model the interactions between the two physics. As an alternative to density methods, level set based methods \cite{osher1988a} show great promise wrt. coupled multi-physics problems. \markline{These methods implicitly defines the geometry by an iso-level of the level set function, which is usually taken as the zero level contour of the level set function.} However, when such methods are used with an ersatz material approach \cite{wang2008a}, the level set function is - similar to density methods -  mapped onto a piece-wise constant density field and the interface is still represented with an interpolated gray area. Hence, such level set methods suffers the same drawbacks as the density based topology optimization for coupled problems. Alternatively, the geometry defined by the zero level of the level set function can also be captured with body fitted meshes \cite{allaire2014a} in which a very accurate modelling of the coupled physics can be realized depending on the quality of the elements along the re-meshed interface. This approach commonly uses the solution of a Hamilton Jacobi type equation to update the design by moving the interface based on the calculated shape sensitivities. However, due to the re-meshing operation at each design iteration, the approach is not efficient for parallel computing frameworks. Moreover, numerical noise may be introduced in the sensitivities due to the extrapolation from the fitted mesh onto the background mesh. Examples of classical level set approach for the optimization of vibroacoustic systems can be found in \cite{shu2014a,isakari2017a,desai2018a}. A detailed comparative review of density, level set and evolutionary based optimization methods for vibroacoustics can be found in the work of \cite{dilgen2019a}.

Cumbersome re-meshing operation can be avoided when the level set approach is coupled with immersed boundary methods \cite{sethian2000a} to capture the geometry. The present work utilizes the immersed boundary cut element method presented in \cite{Schousboe2020} to model the exact (piece-wise linear) boundary represented by the zero contour. 
Throughout the optimization, a fixed mesh is used and the accurate modeling of the interface elements is realized through a special integration scheme. The employed cut element method can be categorized as a special (simplified) case of CutFEM without stabilization \cite{burman2015a}, finite cell method \cite{duester2008a} or X-FEM without enrichment \cite{daux2000a}. For the optimization, the utilized cut element method is coupled with the so-called explicit level set approach \cite{de2004a,kreissl2012a} in which the nodal level set values are directly tied to the mathematical design variables. This approach enables the utilization of nonlinear programming tools, such as the Method of Moving Asymptotes (MMA) algorithm \cite{svanberg1987a} that is used in the current work, and in turn allows for general optimization frameworks where multiple objective and constraint functions easily can be considered in the optimization.

The current work focuses on topology, or generalized shape, optimization of transient vibroacoustic problems where a time dependent problem formulation allows the realization of wideband performance optimization in the frequency domain. Throughout the work, the frequency response of the coupled vibroacoustic system is obtained through a fast Fourier transform (FFT) algorithm applied to the transient response of the system. This allows the objective function to be defined in frequency domain where the optimization can directly tailor the response of the coupled system. To the best of the authors’ knowledge, utilization of time-domain methods to realize wideband optimization in the frequency domain have not been demonstrated before for structural optimization of strongly coupled vibroacoustic problems. Moreover, the work utilizes the discrete adjoint method for calculating the gradients of the objective function \cite{Pollini2017,Dilgen2020_trans_sens}. A thorough explanation of the consistent sensitivity calculation through the discrete adjoint framework that involves the FFT operation, which is needed to define the objective function in frequency domain, is also among the contributions of the present work. \marklineNew{The remainder of the paper is organized as follows.} Section \ref{TheorySection} introduces the level set geometry representation, the vibroacoustic governing equations along with the immersed boundary formulation and spatial and temporal discretizations. 
The optimization problem as well as the sensitivity analysis using the discrete adjoint method are introduced in section \ref{OPTSECTION}. Section \ref{NumSetup} describes the numerical setup that is used for the considered examples throughout the work. The numerical examples are presented in section \ref{NUMERICALEX} where the developed framework is applied for the optimization of various vibroacoustic filter designs.

\section{Theory and methods} \label{TheorySection}
\subsection{Geometry representation}
\label{GeoDescript}
\begin{figure*}
  \input{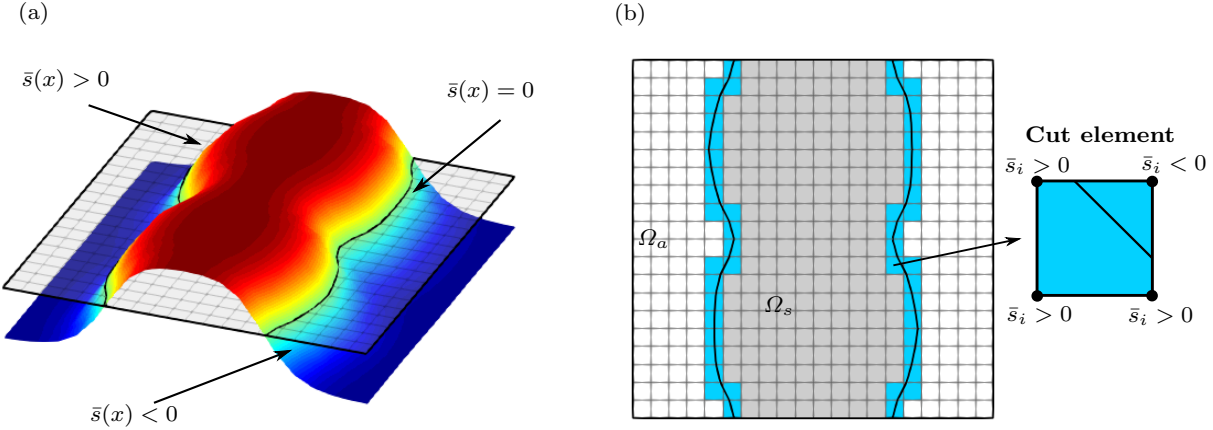}
  \centering
  \caption{An example level set function showing the embedded physical domains. (a) shows the rule that is used for specifying different physical domains embedded with the level set function. (b) finite element mesh where white color shows the uncut elements in the acoustic domain, gray color shows the uncut elements in the solid region and the blue color shows the cut elements.}
  \label{fig:levelsetCut}
\end{figure*}

As the aim of this work is to perform structural optimization for vibroacoustic design problems, it is natural to start by introducing the geometry representation. Throughout this work, the geometry is represented by a scalar valued function $\bar{s}$, i.e. a level set function, which is used to identify the acoustic and structural domains embedded inside the total computational domain $\Omega = \Omega_s \cup \Omega_a$ as illustrated in figure \ref{fig:levelsetCut}(a). This work follows the same geometry representation as used in \cite{Dilgen2020_trans_sens,Dilgen2021}  meaning that the \markline{iso-level is chosen as the zero-level of the level set function and hence, that the}  following rule is used to determine acoustic and structural material phases
\begin{align}
&\bar{s}(x) > 0, \qquad x \in \Omega_s \;\;\text{(structural domain)} \nonumber \\ 
&\bar{s}(x) = 0, \qquad x \in \Gamma_{as} \;\;\text{(interface)} \label{levelDef} \\ 
&\bar{s}(x) < 0, \qquad x \in \Omega_a \;\;\text{(acoustic domain)} \nonumber 
\end{align}

This allows for a versatile geometry description capable of representing complex geometries using the rule given in equation \ref{levelDef} where positive level set values represent the structural domain and negative values identify the acoustic domain. Such an approach is especially well-suited for  acoustic-structure interaction problems, where the zero iso-level of the function $\bar{s}$ provides a clear definition of the interface $\Gamma_{as}$ between the two physical domain.

\subsection{Governing equations} \label{governSection}
\begin{figure*}
  \input{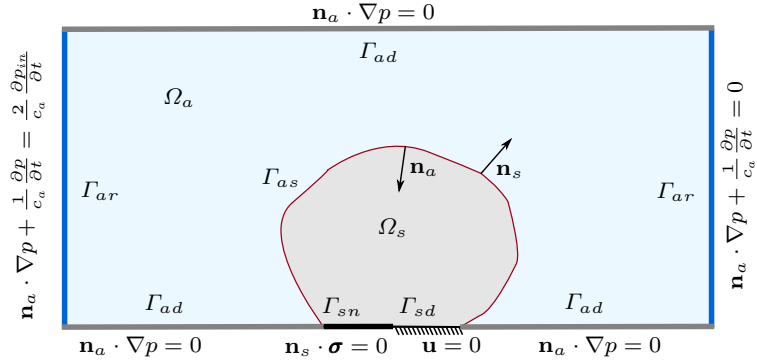}
  \centering
  \caption{Schematic illustration of the coupled acoustic-structural system. Physical domains, the coupled interface and the boundary conditions are showed.}
  \label{fig:physics}
\end{figure*}

The  vibroacoustic system considered in this work consists of plane stress linear elasticity and the linear acoustic wave equation formulated in the time domain as shown in figure \ref{fig:physics}. In order to realize an efficient parallel modelling and optimization framework, we once again emphasize that this work utilizes an immersed boundary cut element method in which the solutions to both physics are obtained in the entire computational domain $\Omega$. This is achieved by application of a fictitious domain approach in which, for the structural domain, a void  phase with low stiffness is defined in the acoustic region $\Omega_a$. For the solution of the acoustic pressure in its fictitious domain (structural domain $\Omega_s$), a rigid phase is defined by modifying the material properties of the acoustic medium.
To this end, considering a transient motion, the linear elasticity equation governing the structural response is written as

\begin{align}
& \rho_s(\mathbf{x}) \frac{\partial^2\mathbf{u}}{\partial t^2} - \nabla \cdot \pmb{\sigma}(\mathbf{x})  +\, \rho_s(\mathbf{x}) \alpha_d \frac{\partial\mathbf{u}}{\partial t} -  \label{st_eq}  \\ 
& \qquad \nabla \cdot \left(\beta_d \frac{\partial \pmb{\sigma}(\mathbf{x}) }{ \partial t} \right)  = 0 \qquad \text{in} \quad \Omega \nonumber \\
& \mathbf{u} = 0  \qquad \qquad \qquad \;\;\;\;\;\; \text{on} \quad \Gamma_{sd} \label{disp_b1}\\
& \mathbf{n}_s \cdot \pmb{\sigma} = 0 \qquad \quad \quad \;\;\;\;\;\;\; \text{on} \quad \Gamma_{sn}\label{disp_b2} \\
& \mathbf{n}_s \cdot \pmb{\sigma} =p\,\mathbf{n}_a  \qquad  \;  \;\;\;\;\;\;\;\;\;  \text{on}  \quad \Gamma_{as}\label{disp_b3}
\end{align}

where $\mathbf{u}$ is the displacement vector, $p$ is the acoustic pressure, $\rho_s(\mathbf{x})$ is the spatially varying density of the solid, $\mathbf{n}_a$ is the normal vector pointing outwards from the acoustic region likewise $\mathbf{n}_s$ is the normal vector defined at the interface pointing outwards from the structural domain. In order to reflect the loss mechanism, the work also considers structural damping using the Rayleigh damping where $\alpha_d$ and $\beta_d$ are the usual Rayleigh damping parameters. Furthermore, the spatially varying Cauchy stress vector $\pmb{\sigma}(\mathbf{x}) $ is defined as $\pmb{\sigma} = \pmb{\mathcal{C}}\left(E_s(\mathbf{x}) ,\nu\right) \pmb{\epsilon}$ in which the constitutive matrix $\pmb{\mathcal{C}}\left(E_s,\nu\right)$ is a function of the spatiallay varying Young's modulus $E_s(\mathbf{x})$ and the Poisson's ratio $\nu$ and the strain vector is defined as $\pmb{\epsilon} = \left[ \frac{\partial u_1}{\partial x}, \;\frac{\partial u_2}{\partial y}, \; \frac{\partial u_1}{\partial y} + \frac{\partial u_2}{\partial x} \right]^\text{T}$. Moreover, equation \ref{disp_b1} defines the fully clamped condition for the structure, equation \ref{disp_b2} specifies the traction free boundary condition whereas the coupling to the acoustic medium is realized through the condition given in equation \ref{disp_b3}. For calculating the transient response of the structure both in structural domain and in its void phase (acoustic domain), the material properties of the solid are altered as
\begin{align}
E_s(\mathbf{x})  = \alpha(\mathbf{x})  \tilde{E}_s, \qquad \rho_s(\mathbf{x})  = \alpha(\mathbf{x})  \tilde{\rho}_s
\end{align}
where the tilde superscript specifies the original material properties and the parameter $\alpha(\mathbf{x}) $ is a dimensionless contrast parameter that is unity for the structural domain and taken as $10^{-8}$ for the void phase of the structure. The acoustic pressure on the other hand is governed by the transient Helmholtz equation which is written as
\begin{align}
& \frac{1}{K_a(\mathbf{x}) } \frac{\partial^2 p}{\partial t^2}  -  \frac{1}{\rho_a(\mathbf{x}) }\nabla^2 p  = 0 \qquad\qquad \text{in} \quad \Omega \\
& \mathbf{n}_a \cdot \nabla p   = 0  \qquad\; \qquad\; \qquad \;\;\; \quad \text{on} \quad \Gamma_{ad}\label{wall}\\
& \mathbf{n}_a \cdot \nabla p   = \rho_a  \frac{\partial^2 \left(\mathbf{n}_{s} \cdot \mathbf{u} \right)}{\partial t^2} \quad  \qquad\; \text{on} \quad \Gamma_{as} \label{coupling}\\
& \mathbf{n}_a \cdot \nabla p + \frac{1}{c_a} \frac{\partial p}{\partial t} = \frac{2}{c_a}\frac{\partial p_{in}}{\partial t} \quad \;\;\;\text{on} \quad \Gamma_{ar} \label{radPbound}
\end{align}

where $\rho_a(\mathbf{x})$ is the spatially varying density of the acoustic medium, $c_a$ is the speed of sound in the acoustic medium. The spatially varying bulk modulus $K_a(\mathbf{x})$  for the acoustic domain is defined as $K_a(\mathbf{x}) =\rho_a(\mathbf{x})  c_a^2$. The boundary conditions used in the current work for the transient acoustic pressure solution are the hard wall condition given in equation \ref{wall}, the coupling boundary condition for the acoustic domain written in equation \ref{coupling} and the absorbing boundary condition given in equation \ref{radPbound}. The absorbing boundary condition is written with a plane wave radiation where $p_{in}$ denotes the transient incoming acoustic wave. Similar to the structural equation, the  transient response of the acoustic pressure is obtained in both the acoustic domain and its rigid phase (structural domain) and is calculated by changing the properties of the acoustic medium as
\begin{align}
K_a(\mathbf{x})  = \frac{\tilde{K}_a}{\alpha(\mathbf{x}) }, \qquad \rho_a(\mathbf{x})  = \frac{\tilde{\rho}_a}{\alpha(\mathbf{x}) }
\end{align}
The dimensionless contrast parameter $\alpha(\mathbf{x})$ is similar to that used for the structural and takes the values of unity for the acoustic domain and $10^{-8}$ for the rigid phase of the acoustic solution. A complete illustration of the resulting coupled acoustic-structure system considered in this work is given in figure \ref{fig:physics} along with all the used boundary conditions. \markline{Note that the use of a single fictitious scaling parameter for both stiffness and density is known to cause issues with spurious modes in the fictitious domains that can lead to unphysical phenomena. However, throughout  extensive numerical experiments, it was found that this did  not pose any problems for the proposed method and we have therefore decided to use a single fixed fictitious domain parameter.}

\subsection{Discretization} \label{discretesection}
This section describes  the spatial and temporal discretizations for the coupled vibroacoustic system introduced in the previous section
%
The spatial discretization is obtained using the standard continuous Galerkin method \cite{zienkiewicz2000a} using Q4 elements and a special integration scheme is employed for the cut elements and the coupling interface, shown in blue and black, respectively, in figure \ref{fig:levelsetCut}. \markline{In short, the utilized cut element approach employs over-integration using a Gaussian quadrature rule for modelling the interface boundary through weighted integration. To this end, a triangulation is carried out in the identified cut Q4 elements to correctly place the integration points in the sub-elements as well as on the interface in order to integrate the coupling boundary conditions, \marklineNew{equations} \ref{disp_b3} and \ref{coupling}, inside the parent Q4 element.} The reader is referred to \cite{Dilgen2020_trans_sens} for details regarding all necessary integrals. The resulting semi-discrete form of the coupled system then reads
\begin{align}
\mathbf{M} \mathbf{\ddot{v}}^n + \mathbf{C} \mathbf{\dot{v}}^n + \mathbf{K} \mathbf{v}^n = \mathbf{h}^n \label{system1}
\end{align}
where $\mathbf{M}$ is the mass matrix, $\mathbf{C}$ is the damping matrix and $\mathbf{K}$ is the stiffness matrix. The superscript $n$ in equation \ref{system1} denotes the current time step. The combined state and load vectors,  $\mathbf{v}$ and $\mathbf{h}$,  are given as

\begin{align}
\mathbf{v} = \begin{bmatrix}
\mathbf{u} \\ \mathbf{p}
\end{bmatrix},\qquad \mathbf{h} = \begin{bmatrix}
0 \\ \mathbf{g}
\end{bmatrix}
\end{align}
where $\mathbf{u}$ and $\mathbf{p}$ are state vectors for structural displacement and acoustic pressure solutions, respectively. The vector $\mathbf{g}$ denotes the source vector for the acoustic equation which is active when the absorption boundary condition (equation \ref{radPbound}) is utilized for radiating a plane wave. The weak form of the coupled system and the identification of the individual element matrices given in equation \ref{system1} are thoroughly described in \cite{Dilgen2020_trans_sens}.

For the temporal discretization, the current work employs an implicit time stepping scheme, namely the Newmark algorithm \cite{Newmark}. The method expresses the first and the second time derivative of the solution vector $\mathbf{v}$ as
\begin{align}
&\mathbf{\dot{v}}^n = a_1 \mathbf{\dot{v}}^{n-1} + a_2 \mathbf{\ddot{v}}^{n-1} + a_3 \left(\mathbf{v}^{n}  - \mathbf{v}^{n-1} \right) \label{1tur} \\
&\mathbf{\ddot{v}}^n = -a_4 \mathbf{\dot{v}}^{n-1} - a_5 \mathbf{\ddot{v}}^{n-1} + a_6 \left(\mathbf{v}^{n}  - \mathbf{v}^{n-1} \label{2tur}\right)
\end{align}
where the parameters $a_1$ to $a_6$ are the so-called Newmark parameters. In order to derive final and complete linear system of equations which is solved in order to determine the solution vector $\mathbf{v}^n$ for the current time $n$, the definitions given in equations \ref{1tur} and \ref{2tur} are substituted into the semi-discrete form of the coupled system (equation \ref{system1}). The resulting equation reads
\begin{align}
&\mathbf{\hat{K}} \mathbf{v}^n = \mathbf{\hat{h}}^n  \label{solveEQ}
\end{align}
where $\mathbf{\hat{K}}$ is an effective stiffness matrix and $\mathbf{\hat{h}}^n$ is the effective load vector, identified as
\begin{align}
&\mathbf{\hat{K}} = \mathbf{K} + a_6 \mathbf{M} + a_3 \mathbf{C} \\
&\mathbf{\hat{h}}^n = \mathbf{h}^n + \mathbf{M} \left( a_4 \mathbf{\dot{v}}^{n-1} + a_5 \mathbf{\ddot{v}}^{n-1} + a_6  \mathbf{v}^{n-1}   \right) + \nonumber \\ &\qquad\mathbf{C} \left( -a_1 \mathbf{\dot{v}}^{n-1} - a_2 \mathbf{\ddot{v}}^{n-1} + a_3 \mathbf{v}^{n-1}   \right) \label{3tur}
\end{align}
Furthermore, the Newmark parameters are identified as
\begin{align}
&a_1 = 1 - \frac{\tilde{\gamma}}{\tilde{\beta}},\qquad a_2 = \left( 1 - \frac{\tilde{\gamma}}{2 \tilde{\beta}} \right) \Delta t, \qquad a_3 = \frac{\tilde{\gamma}}{\tilde{\beta}\Delta t} \\
&a_4 = \frac{1}{\tilde{\beta}\Delta t}, \qquad a_5 = \frac{1}{2 \tilde{\beta}} - 1 , \qquad a_6 = \frac{1}{\tilde{\beta} \Delta t^2} \nonumber
\end{align}
where $\Delta t$ is the time step size. The current work uses a specific family of the Newmark algorithm \cite{jacob1994a} in which the parameters  $\tilde{\beta}$ and $\tilde{\gamma}$ are selected so that the employed time discretization becomes unconditionally stable. $\tilde{\beta}$ and $\tilde{\gamma}$ are given as
\begin{align}
\tilde{\beta} = \frac{1}{4},\qquad \tilde{\gamma} = \frac{1}{2}
\end{align}
The presented work assumes the initial conditions $\mathbf{v}^0$ and $\mathbf{\dot{v}}^0$ are zero which is a minor simplication compared to \cite{Dilgen2020_trans_sens}. The initial condition for the second time derivative of the solution vector $\mathbf{\ddot{v}}^0$ is therefore found by setting the  $\mathbf{v}^0$ and $\mathbf{\dot{v}}^0$ to zero in equation \ref{system1} and solving the following equation
\begin{align}
\mathbf{M} \mathbf{\ddot{v}}^0 = \mathbf{h}^0 \label{initACC}
\end{align}

To simplify notation when formulating the optimization problem, the vectors of state variables and residuals are \markline{collected} in two single vectors $\mathbf{U}^n$ and $\mathbf{R}^n$, respectively, which are given as
\begin{align}
\mathbf{R}^n = \begin{bmatrix}
\mathbf{r}_1^n \\ 
\mathbf{r}_2^n \\ 
\mathbf{r}_3^n
\end{bmatrix},\qquad \mathbf{U}^n = \begin{bmatrix}
\mathbf{v}^n \\
\mathbf{\dot{v}}^n  \\
\mathbf{\ddot{v}}^n 
\end{bmatrix} \label{buyukVec}
\end{align}
Here the individual residual vectors $\mathbf{r}_1^n$, $\mathbf{r}_2^n$ and $\mathbf{r}_3^n$ are identified from the discretized coupled system given in the equations \ref{1tur} to \ref{3tur} as
\begin{align}
&\mathbf{r}_1^n  = \left[ \mathbf{K} + a_6 \mathbf{M} + a_3 \mathbf{C} \right] \mathbf{v}^n - \left[ a_6 \mathbf{M} + a_3 \mathbf{C} \right]\mathbf{v}^{n-1} -  \label{res1} \\ &\qquad \left[a_4 \mathbf{M} - a_1 \mathbf{C} \right]  \mathbf{\dot{v}}^{n-1} + \left[ a_2 \mathbf{C} - a_5 \mathbf{M} \right]\mathbf{\ddot{v}}^{n-1} - \mathbf{h}^{n} \nonumber\\
&\mathbf{r}_2^n = \mathbf{\dot{v}}^{n} - a_1 \mathbf{\dot{v}}^{n-1} - a_2 \mathbf{\ddot{v}}^{n-1} - a_3 \left[ \mathbf{v}^n - \mathbf{v}^{n-1} \right] \\
&\mathbf{r}_3^n = \mathbf{\ddot{v}}^{n} + a_4 \mathbf{\dot{v}}^{n-1} + a_5 \mathbf{\ddot{v}}^{n-1} - a_6 \left[ \mathbf{v}^n - \mathbf{v}^{n-1} \right] \label{res2}
\end{align}

The solution of the coupled vibroacoustic system has been implemented in {\CC} using the PETSc library \cite{petsc-web-page,petsc-user-ref,petsc-efficient} for parallel data management. Moreover, the parallel direct solver MUMPS \cite{MUMPS:1,MUMPS:2} is utilized for the solution of equation \ref{solveEQ}.

\section{Optimization formulation} \label{OPTSECTION}

\subsection{Design parameterization} \label{desparameter}
Having established the physical model and how it is capable of handling arbitrary geometries, the next step is to introduce the design parameterization. To achieve this, a mathematical design variable, $s$, is introduced and then tied to the previously defined level set function $\bar{s}$ through the operations presented here. 
The design variable is, similar to the physical design field $\bar{s}$, defined on nodal points in the computational mesh and is bound between zero and one as this allows for easy use of the MMA optimizer, i.e.
\begin{align}
0 \leq s \leq 1 \label{bounds1}
\end{align}
Although this scaling is appropriate for the optimizer used in this work, it can lead to mesh dependent speeds at which the zero iso-level can change in one design update. This can be alleviated by mapping the design variable into a new variable $\tilde{s}$ which has the upper and lower bounds \marklineNew{corresponding to half the element size \cite{Sharma2017} as}
\begin{align}
-0.5 h_e \leq \tilde{s}(s) \leq 0.5 h_e
\end{align}
Lastly, a smoothing filter is applied to \markline{regularize and} stabilize the optimization process and to obtain the physical design variable $\bar{s}$. This is achieved by use of PDE filter \cite{lazarovFilter} with Neumann boundary conditions, i.e.
\begin{align}
\markline{-r^2 \nabla^2 \bar{s_c} + \bar{s_c} = \tilde{s_c}}
\end{align}
where $r$ is the filter radius \markline{and the subscript $c$ denotes the variables that are defined in cell centers since the filter PDE is solved using a finite volume implementation. Therefore two additional smoothing steps are implied whenever the filter equation is applied due to interpolations between nodes of the mesh and cell centers (from $\tilde{s}$ to $\tilde{s_c}$ and from $\bar{s_c}$ to $\bar{s}$). Here it is noted that the reason the filter PDE is solved with a finite volume implementation is that with using a finite volume method, a numerically stable result can be obtained with any filter radius value (it can even be set to zero).} Finally, it should be remarked that no length scale is imposed as this is deemed outside the scope of the current work since the focus is on  broadband frequency control.

\subsection{Optimization problem formulation} \label{OPTPROBLEM}

To control the structural and acoustic behavior over a wide range of frequencies, we propose to cast the optimization problem in frequency space. This requires that the temporal state response is transformed, which is achieved through the fast Fourier transform (FFT) yielding  the vector of state variables in the frequency domain, $\mathbf{U}^m_f$. Here superscript $m$ refers to frequency point and subscript $f$ is used to distinguish the vector from its temporal counterpart. Here it is noted that due to the FFT operation, the state vector in frequency domain $\mathbf{U}^m_f$ consists of complex numbers. 
Formally, the discrete Fourier transform is defined as
\begin{align}
\mathbf{U}_f^m = \sum_{n=0}^{N-1} \mathbf{U}^n e^{- \frac{i 2 \pi}{N}m n}
\end{align}
where $i$ is the imaginary unit defined as $i=\sqrt{-1}$ and $N$ is the total number of time steps. Note that a Hanning window is applied to the temporal signal before passing it through the FFT and that the current work uses the FFTW library \cite{FFTW98}. 
The generic nested optimization problem considered here can now be formulated as
\begin{align}
&\underset{\mathbf{s}}{\min} \qquad \Phi =  \sum_{m=0}^{M} \phi^m \left( \mathbf{U}^m_f (\mathbf{\bar{s}} )  \right) \label{op1}\\
&\text{s.t.} \qquad \; \mathbf{R}^n\left(\mathbf{\bar{s}} , \mathbf{U}^n (\mathbf{\bar{s}} ) \right) = 0,\qquad \text{for}\;\;\;\;n=0,1,\dots,N  \label{resequ1} \\
&\qquad \;\;\;\;\;\;  \psi_i\left( \mathbf{U}^m_f (\mathbf{\bar{s}} )  \right) \leq 0,\qquad \text{for}\;\;\;\;i=0,1,\dots,K \\
&\qquad \;\;\;\;\;\; s_{min} \leq \mathbf{s} \leq s_{max} \label{op2}
\end{align}
where $\Phi$ is the objective function defined in frequency domain and $M$ is the total number of discrete frequencies that is considered in objective function. Moreover, $\psi_i\left( \mathbf{U}^m_f (\mathbf{\bar{s}} )  \right)$ denotes the $K$ additional constraint functions considered for the optimization. 
The vector of design variables is denoted as $\mathbf{s}$ where $s_{min}$ and $s_{max}$ are its lower and upper bounds (equation \ref{bounds1}) which are set as 0 and 1, respectively. The residual equations are required to be fulfilled at every design step, and is based on the vector of physical design variables $\mathbf{\bar{s}}$ as seen from the equation \ref{resequ1}. 

The optimization problem stated in equations \ref{op1}  to \ref{op2} is solved using the standard Method of Moving Asymptotes (MMA) algorithm \cite{svanberg1987a}. Specifically, the present work includes the parallel implementation of MMA from \cite{aage2013a}.

\subsection{Sensitivity analysis} \label{sensitivitySection}
Sensitivity information is required in order \markline{to} utilize a gradient based optimizer for solving the problem in equations \ref{op1}-\ref{op2}, and this \markline{section }provides the necessary details for a general discrete function defined in the frequency domain based on a discrete temporal response, i.e. $\Phi = \Phi\left( \mathbf{U}_f ( \mathbf{U}( \mathbf{\bar{s}} )) \right)$. We remark that a common trend in time dependent sensitivity analysis is to use the \marklineNew{so-called} semi-discrete adjoint approach, e.g. \cite{dahl2008a}, in which the problem is considered discrete in space but continuous in time. However, the semi-discrete approach is not suitable for cases where the objective function is defined discretely in the frequency domain using FFT. This is because the semi-discrete approach derives the adjoint equation assuming the objective function is defined in time domain. This means that using a discrete FFT cannot be incorporated within the semi-discrete approach. Recently the work of \cite{Dilgen2020_trans_sens} also  showed that even for pure transient optimization problems, the semi-discrete method may not produce consistent sensitivities which can result in gradients having wrong signs. Hence, for calculating consistent and exact sensitivities, the work utilizes a fully discrete sensitivity analysis for calculating the gradient of the objective function with respect to the design variable. However, it should be noted that if a continuous Fourier transform was used, the semi-discrete approach could still be used in which both the temporal and the frequency response are considered as continuous variables \cite{Zhou2021}. 

The derivation of the fully discrete sensitivity analysis is initiated with an augmentation of the objective function by the residual equations and a vector of Lagrangian multipliers $\mathbf{\Lambda}$. Moreover, note that summation of multiple frequencies have been left out in the following to allow for a more condensed description, but that extension to multiple frequencies is \marklineNew{straightforward}. The Lagrangian function $\mathcal{L}$ then reads 
\begin{align}
\mathcal{L} = \Phi \left( \mathbf{U}_f ( \mathbf{U}( \mathbf{\bar{s}} ))  \right) + \sum_{n=0}^{N}  \mathbf{\Lambda}^n \mathbf{R}^n\left(\mathbf{\bar{s}} , \mathbf{U}^n (\mathbf{\bar{s}} ) \right) \label{lagrange}
\end{align}
where it should be highlighted that the objective function is written as a function of discrete state variables in the frequency domain $\mathbf{U}_f$ which are an explicit function of the discrete transient state variables $\mathbf{U}$. 
Due to the Newmark algorithm, similar to the state variables, the vector of Lagrangian multipliers also consists of three fields as
\begin{align}
\mathbf{\Lambda}^n = \begin{bmatrix}
\mathbf{\lambda}^n \\
\mathbf{\dot{\lambda}}^n \\
\mathbf{\ddot{\lambda}}^n 
\end{bmatrix}
\end{align}
The derivative of the Lagrangian function with respect to the physical design variable $\mathbf{\bar{s}}$ can now be written as
\begin{align}
&\frac{\text{d} \mathcal{L}}{\text{d} \mathbf{\bar{s}}} = \sum_{n=0}^{N}  \underbrace{\left( \frac{\partial \Phi}{\partial \mathbf{U}_f} \frac{\partial \mathbf{U}_f}{\partial \mathbf{U}} \right)^n}_{\frac{\partial \Phi}{\partial \mathbf{U}}^n}   \frac{\partial \mathbf{U}^n}{\partial \mathbf{\bar{s}}} \quad + \label{lag1}  \\ & \qquad \qquad {\mathbf{\Lambda}^n}^T \left[  \frac{\partial \mathbf{R}^n}{\partial \mathbf{\bar{s}}} +  \frac{\partial \mathbf{R}^n}{\partial \mathbf{U}^n}   \frac{\partial \mathbf{U}^n}{\partial \mathbf{\bar{s}}} \right] \nonumber
\end{align}
The term $\left( \frac{\partial \Phi}{\partial \mathbf{U}_f} \frac{\partial \mathbf{U}_f}{\partial \mathbf{U}} \right)^n$ on the right hand side of the equation \ref{lag1} is equivalent to the partial derivative of the objective function with respect to the transient state variable, i.e. $ \frac{\partial \Phi}{\partial \mathbf{U}}^n$, which needs to contain the chain rule describing the connection between frequency and time domain. In order to calculate $ \frac{\partial \Phi}{\partial \mathbf{U}}^n$ , firstly the term $ \frac{\partial \Phi}{\partial \mathbf{U}_f}$ is calculated. Since $\mathbf{U}_f$ consists of complex numbers defined in frequency domain, the partial derivative $ \frac{\partial \Phi}{\partial \mathbf{U}_f}$ also consists of complex numbers and the derivative is therefore realized as
\begin{align}
 \frac{\partial \Phi}{\partial \mathbf{U}_f} =  \frac{\partial \Phi}{\partial \mathbf{U}_{f,r}} + i  \frac{\partial \Phi}{\partial \mathbf{U}_{f,i}} \label{eq:partialsensfreq}
\end{align}
where the subscripts $r$ and $i$ denote the real and imaginary parts of a complex number. \markline{Since multiple frequencies are included, one should in fact compute equation \ref{eq:partialsensfreq} for each frequency $m$ in the considered range.} The term $\frac{\partial \mathbf{U}_f}{\partial \mathbf{U}}$ links the frequency domain to time domain. Since the applied discrete Fourier transform operation can be seen as a linear transformation, the partial derivative $\frac{\partial \mathbf{U}_f}{\partial \mathbf{U}}$ describes the discrete Fourier transform operation itself. Hence, the partial derivative of the objective function with respect to the transient state variable $ \frac{\partial \Phi}{\partial \mathbf{U}}^n$ is rewritten as
\begin{align}
\frac{\partial \Phi}{\partial \mathbf{U}}^n = \left( \frac{\partial \mathbf{U}_f}{\partial \mathbf{U}}^T  \frac{\partial \Phi}{\partial \mathbf{U}_f}  \right)^n \label{disfft}
\end{align}
where $\frac{\partial \mathbf{U}_f}{\partial \mathbf{U}}^T$ is simply the application of the inverse discrete Fourier transform operation. The order of operations in equation \ref{disfft} can be stated as the following two steps: (1) calculate the partial derivative $\frac{\partial \Phi}{\partial \mathbf{U}_f}$ and (2) apply the inverse discrete Fourier transform operation on  $\frac{\partial \Phi}{\partial \mathbf{U}_f}$ to obtain $ \frac{\partial \Phi}{\partial \mathbf{U}}^n$, now defined in time domain. Remark that the sensitivity expression becomes real  after application of the inverse FFT. \markline{For completeness, this process can be stated formally using the inverse discrete Fourier transform as
\begin{align}
\frac{\partial \Phi}{\partial \mathbf{U}}^n= \frac{1}{N}\sum_{m=0}^{N-1}  \frac{\partial \Phi}{\partial \mathbf{U}_{\!f}}^{\!m}  e^{\frac{i 2 \pi}{N}m n}
\end{align}
which also shows how multiple frequencies are included.}

Having introduced the chain rule describing the link between frequency and time domains, in order to continue the sensitivity analysis, the derivative of the Lagrangian function with respect to the physical design variable (equation \ref{lag1})  is rewritten as
\begin{align}
&\frac{\text{d} \mathcal{L}}{\text{d} \mathbf{\bar{s}}}= \sum_{n=0}^{N}  {\mathbf{\Lambda}^n}^T \frac{\partial \mathbf{R}^n}{\partial \mathbf{\bar{s}}} + \underbrace{\left[ \frac{\partial \Phi}{\partial \mathbf{U}}^n +  {\mathbf{\Lambda}^n}^T \frac{\partial \mathbf{R}^n}{\partial \mathbf{U}^n}   \right]}_{=0} \frac{\partial \mathbf{U}^n}{\partial \mathbf{\bar{s}}} \label{set_zero}
\end{align}
Using the fact that the Lagrangian vector can be freely chosen, the underlined part in the above equation is set to zero to prevent the calculation of the partial derivative $\frac{\partial \mathbf{U}^n}{\partial \mathbf{\bar{s}}}$. This gives rise to the adjoint equation
\begin{align}
\frac{\partial \mathbf{R}}{\partial \mathbf{U}}^T \mathbf{\Lambda} = -  \frac{\partial \phi}{\partial \mathbf{U}} \label{adjoint_eq}
\end{align}
where the superscript $n$ is dropped to highlight that the adjoint problem in equation \ref{adjoint_eq} contains the all time steps considered. The partial derivative of the residual vector with respect to the state variables $\frac{\partial \mathbf{R}}{\partial \mathbf{U}}$ has the following general form for a linear set of residuals
\begin{align}
\frac{\partial \mathbf{R}}{\partial \mathbf{U}} = \begin{bmatrix}
	\mathbf{A}_0^{\mathbf{U^0}} &                             &        &                                     &                                     &                   \\
	\mathbf{B}_1^{\mathbf{U^0}} & \mathbf{A}_1^{\mathbf{U^1}} &        &                                     &                                     &                         \\ 
 	                            & \ddots                      & \ddots &                                     &                                     &                          \\
 	                            &                             & \ddots & \ddots                              &                                     &                          \\
 	                            &                             &        & \mathbf{B}_{N-1}^{\mathbf{U^{N-2}}} & \mathbf{A}_{N-1}^{\mathbf{U^{N-1}}} &                              \\ 
  	                            &                             &        &                                     & \mathbf{B}_N^{\mathbf{U^{N-1}}}     & \mathbf{A}_N^{\mathbf{U^N}} \\
\end{bmatrix} \label{DRDU_AD}
\end{align}
Here, the subscripts show the current time step and the superscripts denote the corresponding state variables. \markline{Sub-matrices in $\frac{\partial \mathbf{R}}{\partial \mathbf{U}}$ (equation \ref{DRDU_AD}) are written in general form which considers a time integration scheme that depends on the current time step and one previous time step. From this fact, the residual function for a time step $n$ can be written as}
\begin{align}
\mathbf{R}^n = \mathbf{A} \, \mathbf{U}^n + \mathbf{B}\, \mathbf{U}^{n-1}
\end{align}
Generally, as it can be seen from the above equation, the sub-matrices $\mathbf{A}$ and $\mathbf{B}$ remain constant for each time step and can be identified from the residual equations given in equations \ref{res1} to \ref{res2}. Explicit definitions of the sub-matrices $\mathbf{A}$ and $\mathbf{B}$ considering the Newmark algorithm can be found in \cite{Dilgen2020_trans_sens}. Moreover, the sub-matrix $\mathbf{A}_0$ is easily identified from the initial conditions.
Due to the transpose operation applied on $\frac{\partial \mathbf{R}}{\partial \mathbf{U}}$ for the adjoint equation given in equation \ref{adjoint_eq}, the solution of the adjoint equation is realized from reverse pseudo time steps which are written as
\begin{align}
\left(\mathbf{A}\right)^T \mathbf{\Lambda}^{N} &= \;- \frac{\partial \phi}{\partial \mathbf{U}}^{N}   \\
\left(\mathbf{A} \right)^T \mathbf{\Lambda}^{N-1} &=\; - \frac{\partial \phi}{\partial \mathbf{U}}^{N-1} - \,\, \left(\mathbf{B}\right)^T \mathbf{\Lambda}^{N} \\
&\;\;\vdots \nonumber \\ 
\left(\mathbf{A}_{0} \right)^T \mathbf{\Lambda}^{0} &=\; - \frac{\partial \phi}{\partial \mathbf{U}}^{0} - \,\, \left(\mathbf{B} \right)^T \mathbf{\Lambda}^{1} 
\end{align}
After the adjoint equation is solved for the Lagrangian variables $\mathbf{\Lambda}^{n}$, the final sensitivity of the objective function is calculated as
\begin{align}
&\frac{\text{d} \Phi}{\text{d} \mathbf{\bar{s}}} = {\mathbf{\Lambda}^0}^T \left[ \frac{\partial \mathbf{A_0}}{\partial \mathbf{\bar{s}}}\mathbf{U}^0 \right] +   \sum_{n=1}^{N} {\mathbf{\Lambda}^n}^T \left[ \frac{\partial \mathbf{A}}{\partial \mathbf{\bar{s}}}\mathbf{U}^n + \frac{\partial \mathbf{B}}{\partial \mathbf{\bar{s}}}\mathbf{U}^{n-1}  \right] \label{sensCALC}
\end{align}
Considering a Newmark time stepping scheme, the explicit definitions of the partial derivatives $\frac{\partial \mathbf{A_0}}{\partial \mathbf{\bar{s}}}$, 
$\frac{\partial \mathbf{A}}{\partial \mathbf{\bar{s}}}$ and $\frac{\partial \mathbf{B}}{\partial \mathbf{\bar{s}}}$ are given in \cite{Dilgen2020_trans_sens}. As it can be seen from the gradient calculation in equation \ref{sensCALC}, state variables from the solution of the transient coupled vibroacoustic system needs to be stored in order to carry out discrete adjoint method for transient optimization problems. Alternatively, a check-pointing scheme could be applied to reduce memory usage at the cost of an additional forward analysis.
It is noted that the gradient calculation given in equation \ref{sensCALC} is only done for the cut elements since the partial derivatives  $\frac{\partial \mathbf{A_0}}{\partial \mathbf{\bar{s}}}$, 
$\frac{\partial \mathbf{A}}{\partial \mathbf{\bar{s}}}$ and $\frac{\partial \mathbf{B}}{\partial \mathbf{\bar{s}}}$ are zero elsewhere.

\marklineNew{The} above sensitivity analysis provides the gradient of the objective function with respect to the physical design variable $\mathbf{\bar{s}}$, whereas the gradient with respect to the mathematical design variable $\mathbf{s}$ is obtained from application of the following chain rule
\begin{align}
\frac{\text{d} \Phi}{\text{d} \mathbf{s}} = \frac{\text{d} \Phi}{\text{d} \mathbf{\bar{s}}} \frac{\partial \mathbf{\bar{s}}}{\partial \mathbf{\bar{s}_c}} \frac{\partial \mathbf{\bar{s}_c}}{\partial \mathbf{\tilde{s}_c}} \frac{\partial \mathbf{\tilde{s}_c}}{\partial \mathbf{\tilde{s}}} \frac{\partial \mathbf{\tilde{s}}}{\partial \mathbf{s}}
\end{align}
As it can be seen from the above chain rule, the partial derivatives describe the link between the physical design variable $\mathbf{\bar{s}}$ and mathematical design variable $\mathbf{s}$. The partial derivatives in the chain rule are applied in reverse order in order to calculate the gradient of the objective function with respect to the design variable $\frac{\text{d} \Phi}{\text{d} \mathbf{s}}$. Here, the partial derivative $\frac{\partial \mathbf{\tilde{s}}}{\partial \mathbf{s}}$ describes changing the bounds on the mathematical design variable, $\frac{\partial \mathbf{\tilde{s}_c}}{\partial \mathbf{\tilde{s}}}$ describes the interpolating $ \mathbf{\tilde{s}}$ from nodes to the element centers, $\frac{\partial \mathbf{\bar{s}_c}}{\partial \mathbf{\tilde{s}_c}}$ is the chain rule regarding the employed PDE filter in which its implementation details are thoroughly given in \cite{lazarovFilter} and $\frac{\partial \mathbf{\bar{s}}}{\partial \mathbf{\bar{s}_c}}$ describes the interpolation operation from element centers to the nodes of the mesh.

Here it is noted that, for each optimization case that is considered for the current work, the calculated sensitivities are checked against a first order finite difference calculation. In the worst case, the difference between the calculated sensitivity and the finite difference computation was well below $0.1 \%$ since the utilized discrete adjoint method always yields exact and consistent sensitivities.

\section{Numerical setup} \label{NumSetup}
Here, the numerical setup used for the optimization of transient vibroacoustic filter designs is presented. The section also introduces the objective function along with the computational domain and material properties that are used throughout the work. \marklineNew{The goal of the optimization is to design a structure inside an acoustic channel which acts as, at a certain frequency range, a wave stopper while also allowing the incoming wave to pass at a different predefined frequency set which will be achieved with the internal reflections and the underlying acoustic-structural interactions.} In order to define the objective function which is used for the optimization of acoustic filter designs for a wide frequency range, a measure of transmission is needed to be defined. To achieve this, after the transient response of the state variables are obtained from the solution of the transient coupled vibroacoustic system, the transmitted acoustic pressure signal is integrated at the outlet of the considered acoustic duct as
\begin{align}
\hat{p}(t) = \int_{\Gamma_{out}} p(t)\, \text{d} \Gamma
\end{align}
the frequency response of the transmitted acoustic pressure is found by the discrete Fourier transform
\begin{align}
\hat{p}(f) = \text{FFT}\left( \hat{p}(t) \right)
\end{align}
and the transmission $S(f)$ is then defined as
\begin{align}
S(f) = \frac{\hat{p}(f)}{\hat{p}_0(f)} \label{transmission}
\end{align}
where the subscript $0$ denotes the transmitted acoustic pressure when there is no structure in the acoustic duct which is calculated the same way as $\hat{p}(f)$. As it can be seen from the above equation \ref{transmission}, when the value of $S(f)$ is unity for a particular frequency a full transmission is achieved. Meaning that the recorded amplitude of the transmitted acoustic pressure for an empty acoustic duct at a certain frequency is equal to that of an acoustic duct with a vibrating structure present in it. Moreover when the value of $S(f)$ is zero for a certain frequency, there is no transmission realized compared to the transmission of the empty duct. In other words, the vibrating structure does not transmit acoustic pressure towards to outlet of the duct. 

In order to reflect the \marklineNew{so-called} pass-band and stop-band regions in the considered frequency span, two objective functions are considered which are written as
\begin{align}
&\Phi_1 \left( \mathbf{U}_f ( \mathbf{U}( \mathbf{\bar{s}} ))  \right) = \sum_{f=n_1}^{n_2} \frac{\left( S(f) - a \right)^2}{a^2},\qquad a=1 \label{aPar}\\
&\Phi_2 \left( \mathbf{U}_f ( \mathbf{U}( \mathbf{\bar{s}} ))  \right) = \sum_{f=n_3}^{n_4} \frac{\left( S(f) - b \right)^2}{b^2} \label{bPar}
\end{align}
where the minimization of the function $\Phi_1$ fits the frequency response of the transmitted acoustic pressure to that of an empty acoustic duct hence realizes a full transmission in a frequency window defined between $n_1$ and $n_2$. The minimization of $\Phi_2$ on the other hand lowers the amplitude of the frequency response of the transmitted acoustic pressure compared to the response of an empty acoustic duct in a frequency range defined between $n_3$ and $n_4$. The order of magnitude difference between the transmitted acoustic pressure in the stop-band region and the empty acoustic duct is defined by the parameter $b$ in equation \ref{bPar}. The study on the selection of $b$ in order to realize an effective zero-transmission and its effect on the overall optimization performance are given in section \ref{objStudy}. \markline{Remark that the proposed functions in equations \ref{aPar} and \ref{bPar} were chosen from several candidate functions as they performed best for the filter design problem. \marklineNew{Among the tested candidate formulations, we also tried to minimize the unweighted difference between $S(f)$ and  $a$ together with the difference between $S(f)$ and $b$ and an even simpler approach in which $S(f)$ was maximized for the pass bands and minimized for the stop bands. However, none of these approaches worked as well as the weighted least squares type formulation.} Moreover, it should be emphasized that solving dynamic problems over wide frequency ranges presents a truly challenging design problem with many local minima due to the non-convexity. Hence, there may be better problem formulations available, but out of all the formulation investigated, the one presented next is the one that proved best for the filter design problems.}

In order to realize an optimization problem where multiple objective functions are equally minimized, the problem can be formally written in a so-called min-max formulation as
\markline{
 \begin{align}
 &\underset{\mathbf{s}}{\min} \qquad \max\{\Phi_1 \left( \mathbf{U}_f ( \mathbf{U}( \mathbf{\bar{s}} ))  \right),\Phi_2 \left( \mathbf{U}_f ( \mathbf{U}( \mathbf{\bar{s}} ))  \right)\}  \\
 &\text{s.t.} \qquad \; \mathbf{R}^n\left(\mathbf{\bar{s}} , \mathbf{U}^n (\mathbf{\bar{s}} ) \right) = 0,\qquad \text{for}\;\;\;\;n=0,1,\dots,N  \nonumber \\
  &\qquad \;\;\;\;\;\; 0 \leq \mathbf{s} \leq 1 \nonumber 
 \end{align}}
However, since the min-max formulation is not differentiable an additional variable is introduced and the optimization problem is cast as a bound formulation as
 \begin{align}
 &\underset{\mathbf{s},\,z}{\min} \qquad z \label{boundE1}\\
 &\text{s.t.} \qquad \; \mathbf{R}^n\left(\mathbf{\bar{s}} , \mathbf{U}^n (\mathbf{\bar{s}} ) \right) = 0,\quad \text{for}\;n=0,1,\dots,N   \\
&\qquad \;\;\;\;\;\; \Phi_1 \left( \mathbf{U}_f ( \mathbf{U}( \mathbf{\bar{s}} ))  \right) < z \\
 &\qquad \;\;\;\;\;\; \Phi_2 \left( \mathbf{U}_f ( \mathbf{U}( \mathbf{\bar{s}} ))  \right) < z \\
 &\qquad \;\;\;\;\;\; 0 \leq \mathbf{s} \leq 1 \label{boundE2}
\end{align}
where the additional variable \markline{z} is the upper bound for the optimization. The formulation realizes the minimization of the upper bound \markline{$z$} where the objective functions $\Phi_1$ and $\Phi_2$ are defined as constraints in the optimization problem, hence effectively minimizing both $\Phi_1$ and $\Phi_2$. \markline{In what follows, the functions $\Phi_1$ and $\Phi_2$ are therefore referred to as constraint functions.  Note that no external box constraints are imposed on $z$ as this is included using the internal bound variable in MMA.}

\begin{figure*}
  \input{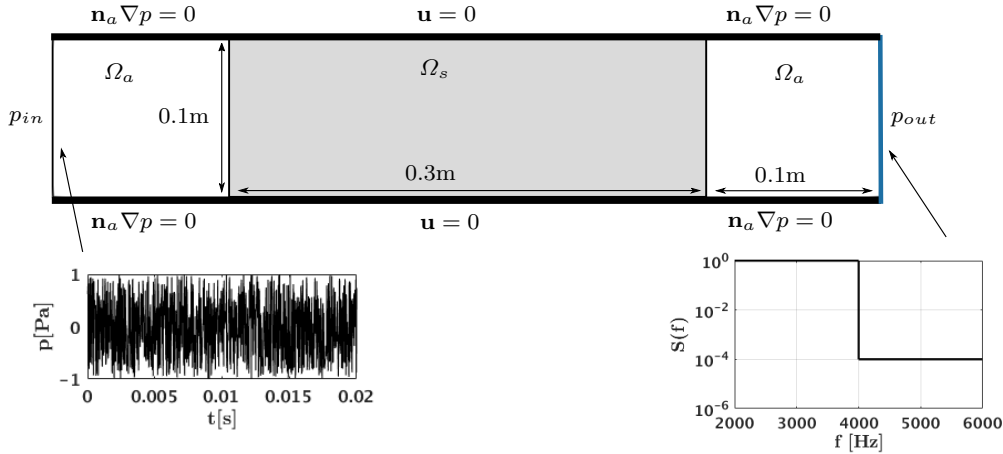}
  \centering
  \caption{Schematic illustration of the optimization case for the design of acoustic filters showing the boundary conditions of the optimization problem. Gray color shows the design domain. Incoming acoustic white noise and an example desired filter shape at the outlet are also shown.}
  \label{fig:desDes}
\end{figure*}

A schematic illustration of the numerical setup for the optimization of acoustic filter designs is given in figure \ref{fig:desDes}. As it can be seen from the figure, both top and bottom boundaries of the acoustic duct are set as hard-wall condition for the acoustic pressure while the structure is considered to be clamped. For the acoustic domain, both right and \marklineNew{leftmost} boundaries are treated as absorbing conditions. The incoming acoustic plane wave from the \marklineNew{leftmost} boundary is also shown in the figure. In order to excite broad frequency range in the coupled system, the incoming plane wave is realized as a white noise with random acoustic pressure values between $-1\,\rm{Pa}$ and $1\,\rm{Pa}$. Moreover, the figure also illustrates an example filter shape at the \marklineNew{rightmost} boundary (at the acoustic duct outlet) where the acoustic transmission of the design $S(f)$ is calculated in frequency domain and a certain filter shape is applied to it. Here it is noted that the presented filter shape in the schematic illustration is given as an example. The work considers various different filters for the optimization which are presented in the following sections.

\begin{figure}
        \centering
  \includegraphics[width=0.5\textwidth]{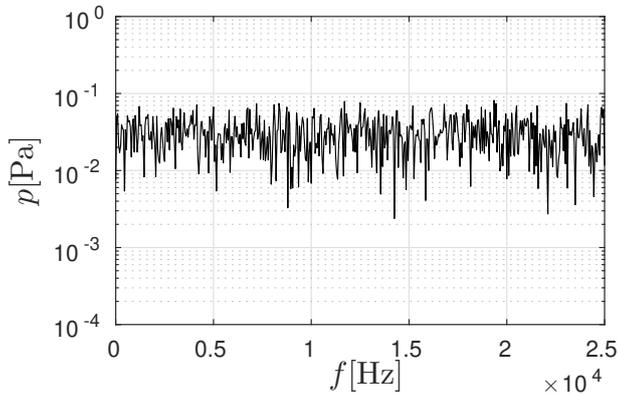}
  \caption{Frequency content of the incoming white noise for the optimization of acoustic filters.}
  \label{fig:example_Freq}
\end{figure}

The overall calculation time for the considered cases throughout the work is chosen as $0.02\rm{s}$ as it can also be seen from the incoming transient acoustic plane wave plot given in figure \ref{fig:desDes}. In order to adequately resolve the coupled vibroacoustic problem in time, the time step size is set to $ \Delta t = 2 \times 10^{-5}\rm{s}$. The frequency content of the white noise applied at the inlet of the acoustic duct to excite broad frequency range in the system is given in figure \ref{fig:example_Freq}. As it can be seen from the figure, the incoming acoustic wave has an approximately constant energy content across all the frequencies present in the signal.

For the optimization, the acoustic domain is taken as air while the structure is considered to be a rubber-like material. The material properties used for the structure and acoustic domains are listed in the tables \ref{tab:MatParam3} and \ref{tab:MatParam4}, respectively.

\begin{table}
        \centering
    \begin{tabular}{ccc}
      \toprule
        $E\;[\rm{Pa}]$ & $\nu $ & $\rho_s \;[\rm{kg}/\rm{m}^3]$   \\
      \midrule
        $50 \times 10^6$ & $0.4$ & $1000$    \\
      \bottomrule
    \end{tabular}
     \caption{Material properties considered for the structure.}
  \label{tab:MatParam3}
\end{table}

\begin{table}
        \centering
    \begin{tabular}{cc}
      \toprule
        $c_a\;[\rm{m/s}]$ & $\rho_a \;[\rm{kg}/\rm{m}^3]$   \\
      \midrule
        $343$ & $1.21$   \\
      \bottomrule
    \end{tabular}
     \caption{Material properties considered for the acoustic domain.}
  \label{tab:MatParam4}
\end{table}

As described in section \ref{governSection}, the work also considers Rayleigh damping for structural damping. Damping for structure is utilized both for the added stabilization of modelling and the subsequent optimization and to reflect the loss mechanism of the real world. The Rayleigh parameters $\alpha_d$ and $\beta_d$ are calculated according to \cite{puthanpurayil2016a} as
\begin{align}
\alpha_d = 2 \zeta \frac{\omega_1 \omega_2}{\omega_1 + \omega_2} \label{alpD} \\
\beta_d = 2 \zeta \frac{1}{\omega_1 + \omega_2} \label{betD}
\end{align}  
here $\zeta$ is called the damping ratio and is taken to be $\zeta=0.1$. Moreover, the work assumes that the two natural frequencies $\omega_1$ and $\omega_2$ are $\omega_1 = 1600 \,2 \pi \, \rm{rad/s}$ and  $\omega_2 = 2200 \,2 \pi \, \rm{rad/s}$. Throughout the work, the computational domain given in figure \ref{fig:desDes} is meshed with \marklineNew{structured quad elements in which each element has an edge length} of $2\times 10^{-3} \rm{m}$\markline{, i.e. a total of 12500 elements}. For optimization, the filter radius $r$ is set to $8 \times 10^{-3}\rm{m}$. 

Throughout the work, the utilized MMA algorithm for solving the optimization problem presented in equations \ref{boundE1} to \ref{boundE2} uses the asymptote parameters of $0.5$, $0.7$ and $1.2$ which are used for controlling the initial adaptation, decrease and increase of the asymptotes, respectively. The penalty parameter which is used for constraints in MMA algorithm is chosen as $1000$. For all of the optimization cases that are considered, the work does not consider a specific stopping criteria. \markline{Instead, the optimization is simply run for a fixed number of iterations. This is simply choice is made due to the following reasons. First, since the performance of the design is extremely sensitive to the small design perturbations, as is the case for many dynamic optimization problems \cite{christiansen2016a}, small oscillations in the constraint functions occur throughout the optimization process. Secondly, as we use the MMA $z$ for the bound variable, we do not have explicit control on the move limits for the objective variable, which can also lead to oscillatory behavior. Hence,  the presented optimization results are obtained using a fixed number of iterations. } 
Moreover, the considered optimization problem does not include a volume constraint on the structure. Since, having only a slab of solid material in the design domain does not form a trivial answer to the optimization problem.

\section{Numerical examples} \label{NUMERICALEX}

To demonstrate the capabilities and limitations of the proposed optimization framework, a series of numerical experiments are conducted in order to examine the performance of the chosen \markline{constraint} function, the dependence on initial configurations, its application to filter design as well as a validation study using commercial finite software. \markline{Remark that the final validation of the proposed time-domain cut element optimization framework and its optimization result is obtained by comparing the response to the result of a COMSOL analysis using a fully coupled vibroacoustic time harmonic (steady state) simulation.}

\subsection{\markline{Constraint} function study} \label{objStudy}
The considered \markline{constraint} function, i.e. equations \ref{aPar} and \ref{bPar}, is first studied in order to illuminate the influence of the input parameters.
As it is introduced in the previous section \ref{NumSetup}, the pass-band and stop-band regions for the \markline{constraint} functions are controlled with the parameters $a$ and $b$, respectively. Minimizing the \markline{constraint} function $\Phi_1$ with having the $a$ parameter as unity fits the calculated transmission value $S(f)$ of the structure inside an acoustic duct to $1$, effectively realizing full transmission in the considered frequency range. Ideally, the parameter $b$ on the other hand needs to be set as $0$ in order to realize zero transmission $S(f)=0$ in the frequency range defining the stop-band region of the filter that is considered for the optimization. However, as it can be seen from the equations \ref{aPar} and \ref{bPar}, an inverse weighting is utilized in the \markline{constraint} functions. \marklineNew{Throughout} our numerical experiments it has been found that the utilized inverse weighting provided designs with superior performances compared to the designs obtained using \markline{constraint} functions without any weighting. Hence, a small positive number is used for the $b$ parameter in order to avoid division by zero in the \markline{constraint} function $\Phi_2$. This section investigates the effect of the $b$ parameter on the optimization and the final performance of the optimized design.
\begin{figure*}
\centering
  \subfloat[{} \label{fig:ABinit:a}] {\includegraphics[width=0.8\textwidth]{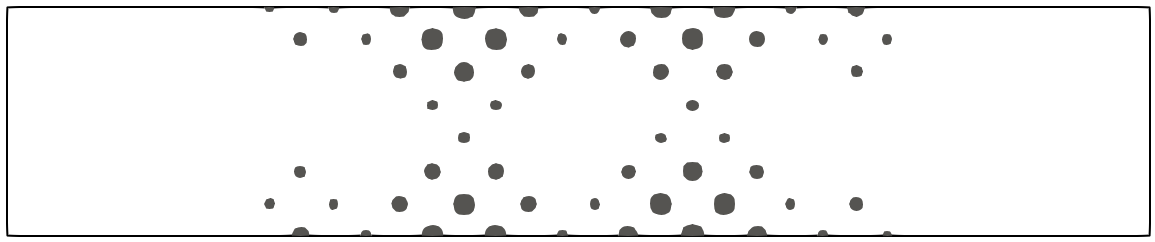}}\\
   \subfloat[{} \label{fig:ABinit:b}] {\includegraphics[width=0.43\textwidth]{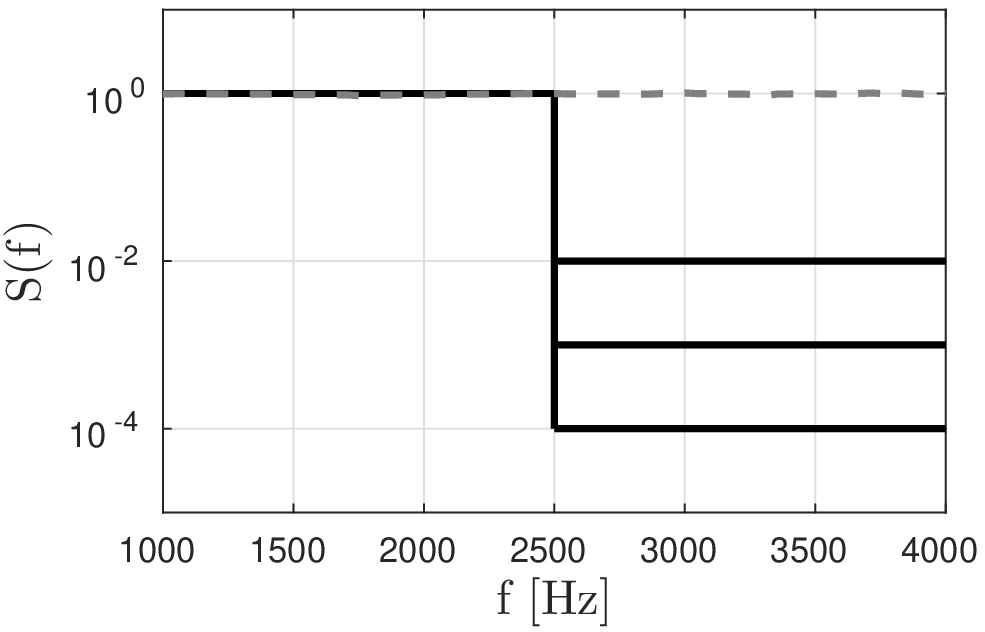}}
   \caption{Optimization for an acoustic low-pass filter design. (a) initial configuration used in the study. (b) Figure showing the low-pass filters for the study where three different levels of transmission $S(f)$ values are considered for the stop-band of the filters which are $1\times 10^{-2}$, $1\times 10^{-3}$ and $1\times 10^{-4}$ (shown with black solid lines). Figure also shows the transmission $S(f)$ of the initial configuration with a dashed gray line.}
  \label{fig:ABinit}
\end{figure*}

Moreover, the section considers a low-pass acoustic filter design for the study. The \markline{constraint} function $\Phi_1$ defining the pass-band region operates on the frequencies between $1000\;\rm{Hz}\leq f \leq 2500\,\rm{Hz}$. For the stop-band, the \markline{constraint} function $\Phi_2$ on the other hand is chosen to be active on the frequencies between $2500\;\rm{Hz} < f \leq 4000\,\rm{Hz}$. Three different $b$ parameters are selected to carry out the study which are $1\times 10^{-2}$, $1\times 10^{-3}$ and $1\times 10^{-4}$. \markline{Remark that several additional values of $b$} \marklineNew{were} \markline{considered in our numerical experiments, but were discarded as they led to poor numerical performance. Further discussion of this follows at the end of this section.}

\begin{figure*}
\centering
  \subfloat[{$\Phi_1=5.06428,\Phi_2= 4.99466$} \label{fig:AB:a}] {\includegraphics[width=0.8\textwidth]{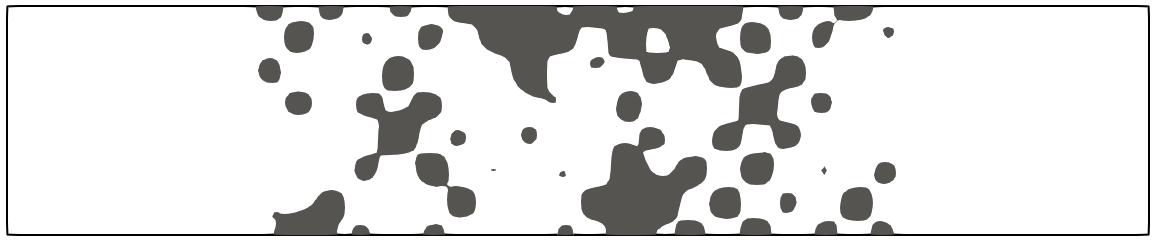}}\\
    \subfloat[{$\Phi_1=10.3386,\Phi_2= 10.5572$} \label{fig:AB:b}] {\includegraphics[width=0.8\textwidth]{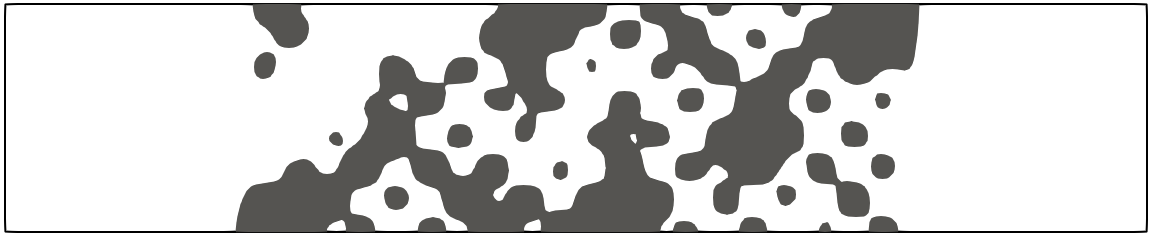}}\\
      \subfloat[{$\Phi_1=14.5365,\Phi_2= 16.9208$} \label{fig:AB:c}] {\includegraphics[width=0.8\textwidth]{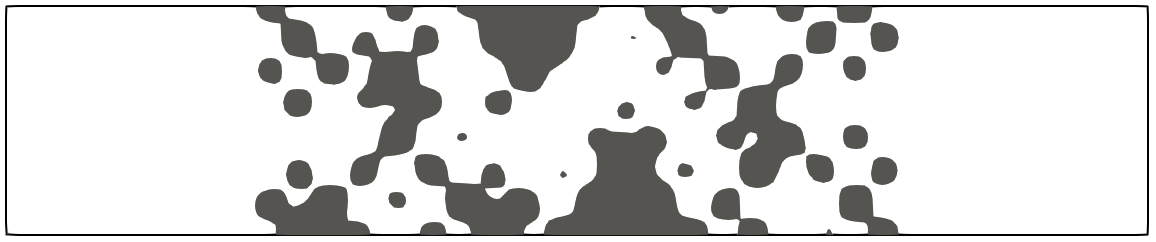}}\\
            \subfloat[{} \label{fig:AB:d}] {\includegraphics[width=0.43\textwidth]{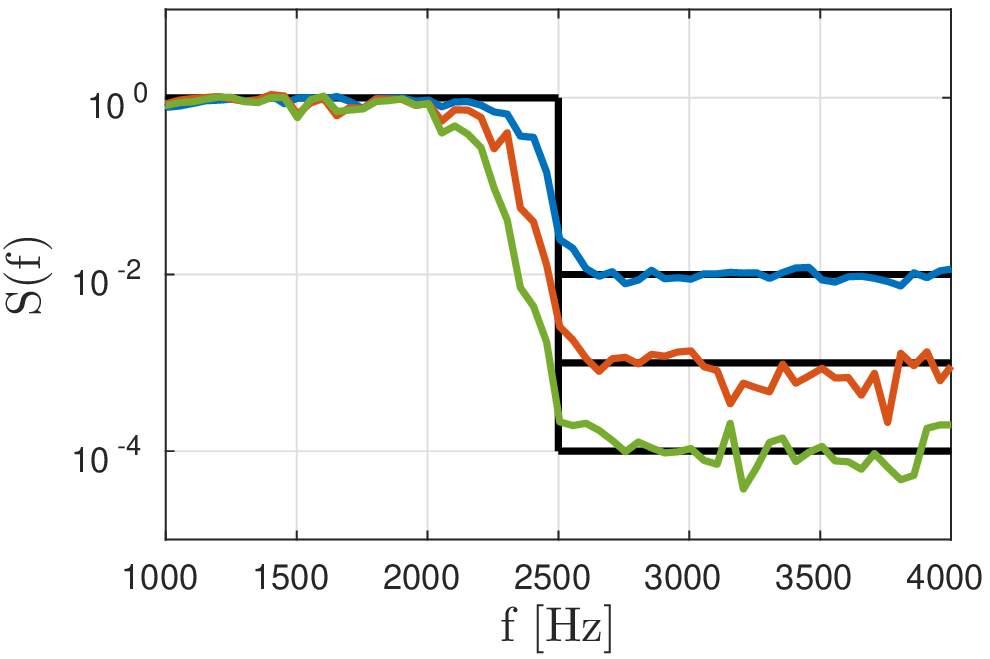}}
                  \subfloat[{} \label{fig:AB:e}] {\includegraphics[width=0.43\textwidth]{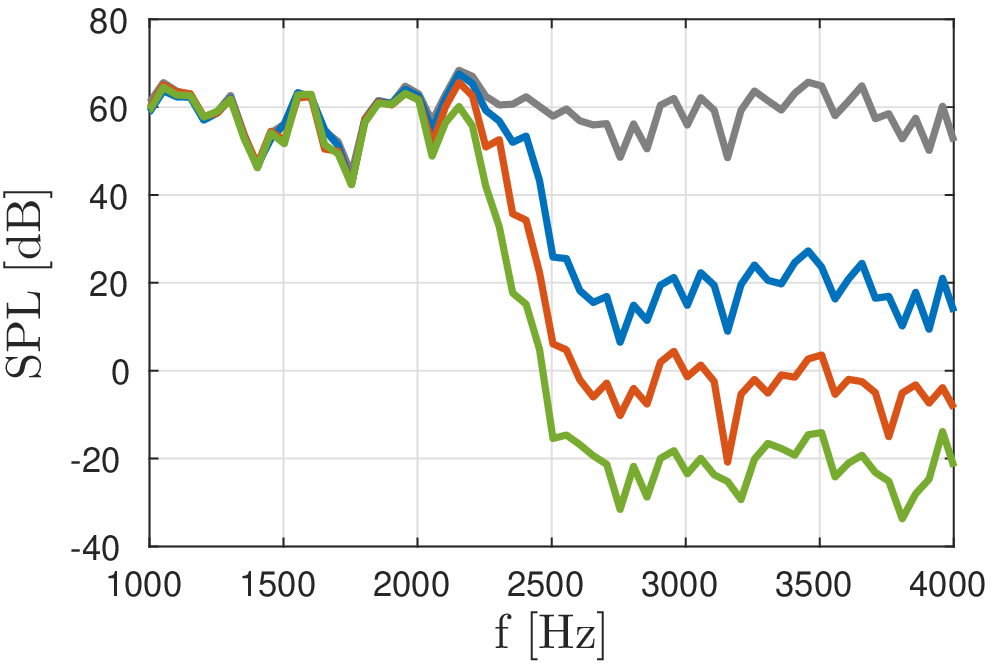}}
   \caption{Optimization for an acoustic low-pass filter design. (a) optimized design for $b=1\times 10^{-2}$. (b) optimized design for $b=1\times 10^{-3}$. (c) optimized design for $b=1\times 10^{-4}$. (d) transmission of the optimized designs in the considered frequency range. Black line is the desired low-pass filter, blue line is the response of the design shown in Fig. \ref{fig:AB:a}, orange line is the response of the design shown in Fig. \ref{fig:AB:b}, green line is the response of the design shown in Fig. \ref{fig:AB:c}. (e) averaged SPL response of the designs calculated at the output, gray line is the empty acoustic duct, blue line is the response of the design shown in Fig. \ref{fig:AB:a}, orange line is the response of the design shown in Fig. \ref{fig:AB:b}, green line is the response of the design shown in Fig. \ref{fig:AB:c}.}
  \label{fig:AB}
\end{figure*}

\begin{figure*}
\centering
 \subfloat{\includegraphics[width=0.4\textwidth]{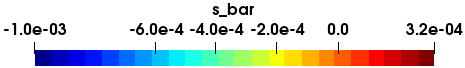}} \addtocounter{subfigure}{-1}  \\
  \subfloat[{} \label{fig:CONGRAPH:a}] {\includegraphics[width=0.8\textwidth]{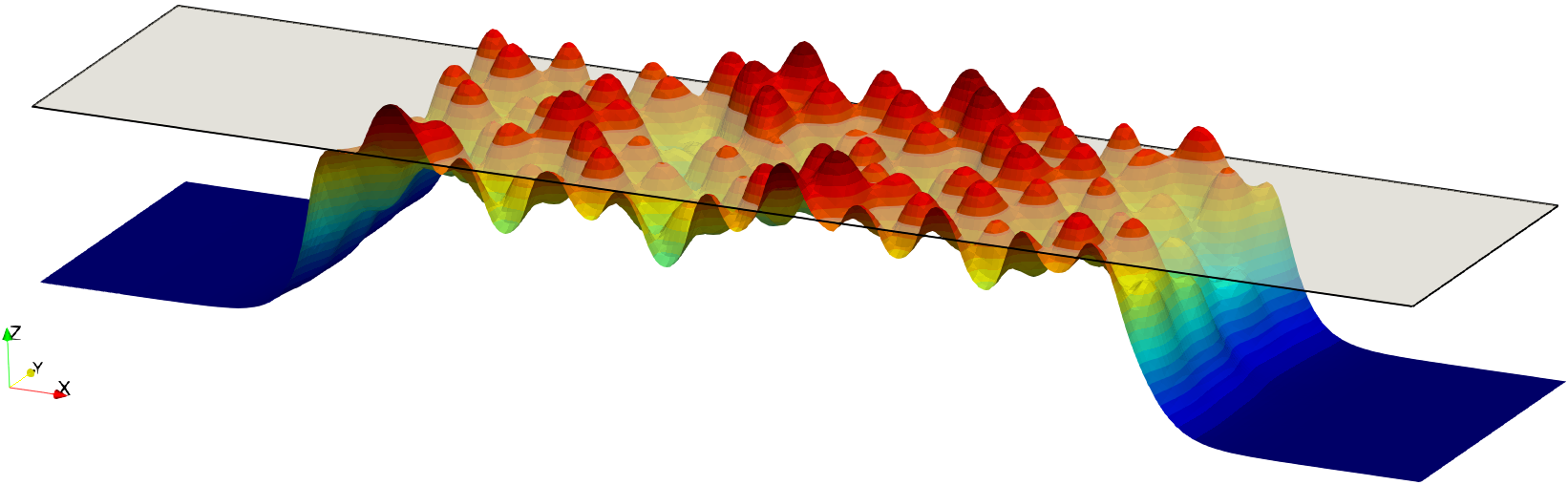}}\\
   \subfloat[{} \label{fig:CONGRAPH:b}] {\includegraphics[width=0.5\textwidth]{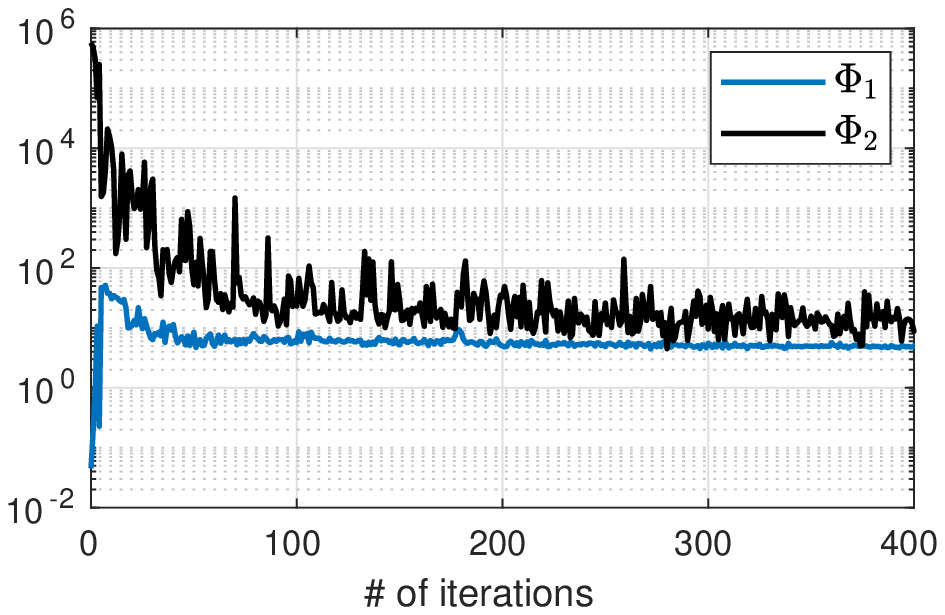}}
   \caption{(a) Optimized level set function $\bar{s}$ of the design given in figure \ref{fig:AB:a}, also showing the zero iso-level of the level set function from which the optimized design is obtained. (b) Constraint function history for the optimized design given in figure \ref{fig:AB:a} run for 400 iterations.}
  \label{fig:CONGRAPH}
\end{figure*}

\begin{table}
        \centering
    \begin{tabular}{ccc}
      \toprule
        $b$ & $\Phi_1 $ & $\Phi_2$   \\
      \midrule
        $1 \times 10^{-2}$ & $0.0455159$ & $576902$    \\
        $1 \times 10^{-3}$ & $0.0455159$ & $5.8754\times 10^7$    \\
        $1 \times 10^{-4}$ & $0.0455159$ & $5.8860 \times 10^9$    \\        
      \bottomrule
    \end{tabular}
     \caption{Calculated \markline{constraint} values $\Phi_1$ and $\Phi_2$ for various $b$ parameters.}
  \label{tab:MatParam6}
\end{table}

Figure \ref{fig:ABinit:a} shows the initial configuration that is used for the optimization. The physical design variables $\mathbf{\bar{s}}$ that specifies the initial design is obtained first by calculating the following expression
\begin{align}
\mathbf{s}_v = \cos \left( \frac{r_1 \pi \mathbf{x} }{l_x} \right) \,\cos \left( \frac{r_2 \pi \mathbf{y} }{l_y} \right) + 0.1 \label{initialEquationS}
\end{align}
where $\mathbf{x}$ and $\mathbf{y}$ are the nodal coordinates of the mesh in the design domain, $l_x$ and $l_y$ is taken to be $0.1$ and $r_1$ and $r_2$ is chosen as $7$. After calculating $\mathbf{s}_v$, the mathematical design variables $\mathbf{s}$ are obtained with the following rule as

\begin{align}
s_i = \begin{cases}
    0 ,& \text{if}\qquad s_{v,i} \geq 0.01\\
    1 ,& \text{if}\qquad s_{v,i} < 0.01
\end{cases}
\end{align}
With applying the design parameterization described in section \ref{desparameter}, the physical design variables $\mathbf{\bar{s}}$ containing the initial design is obtained. Physically, the initial design given in figure \ref{fig:ABinit:a} can be considered as an infinitely long acoustic duct in the out of plane direction where the two dimensional approximation realizes its cross section view from the mid-section of the channel.

Figure \ref{fig:ABinit:b} shows the desired acoustic low-pass filters with the three different levels of stop-band in which the transmission $S(f)$ is to be fitted to the values of $1\times 10^{-2}$, $1\times 10^{-3}$ and $1\times 10^{-4}$. As it can be seen from the figure, the initial configuration nearly has a full transmission across the considered frequency range of $1000\,\rm{Hz}$ to $4000\,\rm{Hz}$. Table \ref{tab:MatParam6} lists the calculated \markline{constraint} values for the initial configuration for each optimization case. As expected, since the same initial design is utilized and the parameter $a$ is set to $a=1$ for each case, the values of $\Phi_1$ are the same. However, the values of $\Phi_2$ increase two orders of magnitude for each order of magnitude decrease in the $b$ parameter. Since, the \markline{constraint} values $\Phi_1$ and $\Phi_2$ are inversely weighted with the square of $a$ and $b$, respectively.

\begin{figure*}
\centering
  \subfloat[{} \label{fig:ABFail:a}] {\includegraphics[width=0.8\textwidth]{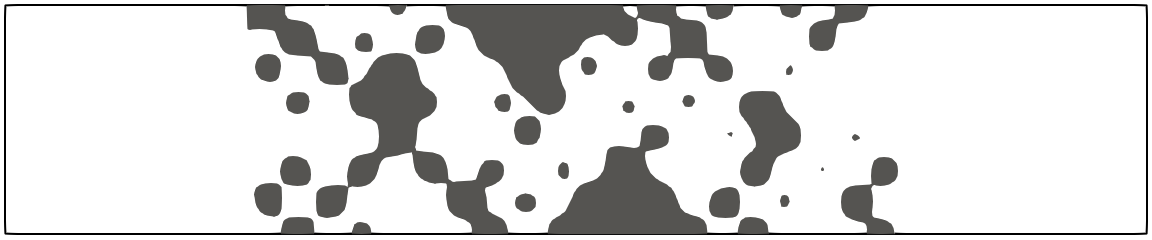}}\\
   \subfloat[{} \label{fig:ABFail:b}] {\includegraphics[width=0.43\textwidth]{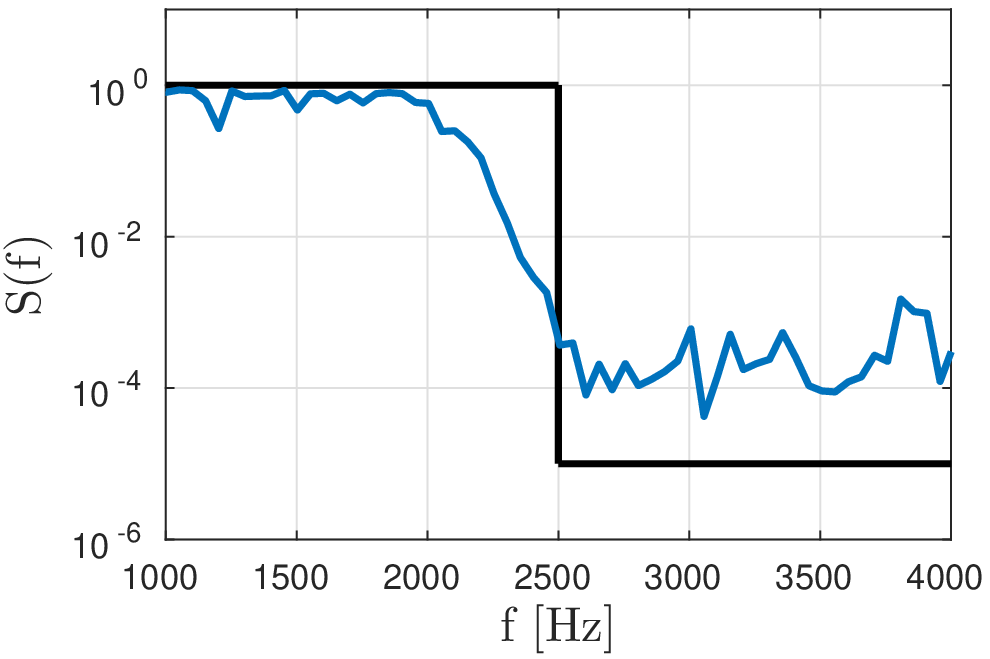}}
   \caption{Optimization for an acoustic low-pass filter design. (a) optimized design for $b=1\times 10^{-5}$. (b) transmission of the optimized design in the considered frequency range. Black line is the desired low-pass filter, blue line is the response of the optimized design.}
  \label{fig:ABFail}
\end{figure*}

The result of the comparative study is given the figure \ref{fig:AB}. As it can be seen from the figures \ref{fig:AB:a} to \ref{fig:AB:c}, although the optimized designs have highly complex topologies, a similar trend can be seen in how the structures clustered inside the design domain for each case. Also, due to the lack of feature size control in the optimization, minimum feature size is bound by the element size used in the computational mesh which can be seen in the thin connections and small island structures present in the optimized designs. 
\markline{Interestingly, the optimized low pass filter designs share little resemblance to classical low-pass filters, in which the channel width is extended for a short duration. Based on our numerical experiments, we conjecture that the performance of the optimized design relies on several, individual resonators and reflectors. However, visualizing this phenomena using the transient, multi-frequency input signal is quite difficult, if not impossible, and the reader is referred to the post-validation in section \ref{sec:validation} where a low pass filter is analyzed in COMSOL using a time-harmonic analysis.}

The frequency responses of the optimized designs can be seen in the figure \ref{fig:AB:d} where the calculated transmission $S(f)$ at the outlet is plotted for all optimized designs. It can be seen from the figure that the the calculated transmission $S(f)$ for each case respects the desired low-pass filter shape for the considered frequency range. However, as the parameter $b$ is decreased for having a lower transmitted acoustic pressure response at the stop-band region, the resulted optimized designs' performances for the transition between the pass-band and stop-band regions also decrease. This effect can also be seen in the presented \markline{constraint} function values $\Phi_1$ and $\Phi_2$ for each optimized design in figures \ref{fig:AB:a} to \ref{fig:AB:c} where the \markline{constraint} values increased as the transmission value is decreased in the stop-band of the acoustic filter. As it can also be seen from the presented \markline{constraint} values, optimized designs' $\Phi_1$ values actually end up at higher values compared to the initial configuration which is due to the transition between the pass-band and stop-band regions for the considered low-pass filter.

Overall, compared to the initial transmission given in figure \ref{fig:ABinit:b}, it can be seen from the figure \ref{fig:AB:d} that the developed optimization framework can successfully tailor the frequency content obtained from the transient response of the coupled vibroacoustic system, effectively designing acoustic filters. Moreover, figure \ref{fig:AB:e} compares the averaged sound pressure level (SPL) values at the outlet of the acoustic duct for the considered frequency range. Compared to the SPL response of the empty acoustic duct, the considered low-pass filter shape can also be seen from the SPL values of each optimized design. Here, the optimized designs have nearly the same averaged SPL response with the empty acoustic duct from $1000\,\rm{Hz}$ to approximately $2250\,\rm{Hz}$ and after a sharp transition into the stop-band region, transmission $S(f)$ values of $1\times 10^{-2}$, $1\times 10^{-3}$ and $1\times 10^{-4}$ roughly corresponds to the averaged SPL values of $20\,\rm{dB}$, $0\,\rm{dB}$ and $-20\,\rm{dB}$ for the stop-band region of the low-pass acoustic filter, respectively. 
\markline{To demonstrate what happens if choosing too small a value of $b$,} the same optimization setup is again considered with a low-pass filter design where the $b$ parameter is lowered to $b=1\times 10^{-5}$. Figure \ref{fig:ABFail:a} shows the optimized design where a similar trend is seen compared to the design presented in figure \ref{fig:AB:c}. However, as it can be seen from the performance of the acoustic filter given in figure \ref{fig:ABFail:b}, optimization resulted in a low quality local-minimum. Meaning that the optimized design fails to lower the transmission to $S(f)=1\times 10^{-5}$ for the stop-band region of the considered low-pass filter. \markline{If, on the other hand, $b$ is increased the drop in SPL is too small to yield a functioning filter. Thus,} considering the SPL values given in figure \ref{fig:AB:e} and the corresponding performances for each $b$ parameter, $b=1\times 10^{-3}$ is deemed the most effective for optimizing acoustic filters within the developed framework and will be used for the remaining of the work.

\markline{Finally, in order to illustrate the evolution of the constraint functions' ($\Phi_1$ and $\Phi_2$) throughout optimization iterations, figure \ref{fig:CONGRAPH:b}, shows the constraint function history for the optimized design given in figure \ref{fig:AB:a}. As it can be seen from the figure, the optimization quickly gets the $\Phi_1$ to a stable level, which governs the performance in the pass-band, and realizes a significant reduction in the constraint function $\Phi_2$ and thus achieving a good performance in the stop-band. It is important to stress, that although the constraint function history does contain noticeable oscillations, the design itself converges to the final topology within approximately 300 iterations. The oscillations are therefore a consequence of the fact that even small changes in design can lead to significant changes in the constraint functions
, which as already stated is a known issue with dynamics and structural optimization. To show that the constraint oscillations are not a consequence of an ill-defined level set function, i.e. a  flat level set function for which even very small design variable changes will lead to large changes in the physical design, figure  \ref{fig:CONGRAPH:a} show the final level set function. From this figure it is observed that the optimized level set function is free of regions with flat slopes which shows that the employed filtering is sufficient for the purpose of the presented work. Note that the flat regions at the inlet and outlet are by construction as these regions are non-designable.} 
\markline{Finally, remark that all of the presented optimization cases have very similar objective/constraint function behaviours and for the sake of brevity, the constraint function history graph is only presented for one example.}

\subsection{Initial guess study}
\begin{figure*}
\centering
  \subfloat[{$\Phi_1=0.0154458,\Phi_2= 5.71289\times 10^{7}$} \label{fig:InitS:a}] {\includegraphics[width=0.8\textwidth]{ABinit}}\\
    \subfloat[{$\Phi_1=1.0096,\Phi_2= 5.89415\times 10^{7}$} \label{fig:InitS:b}] {\includegraphics[width=0.8\textwidth]{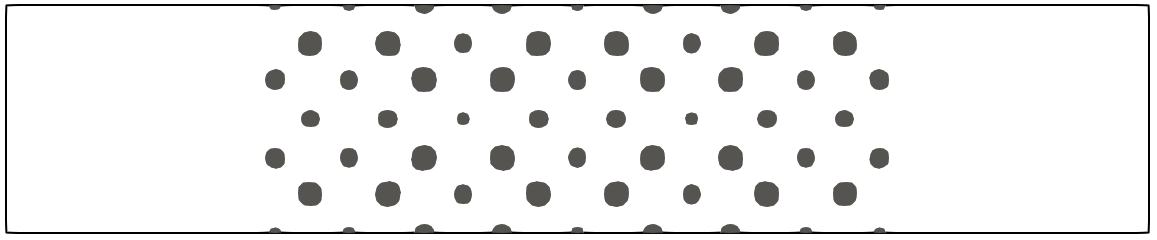}}\\
    \subfloat[{$\Phi_1=0.0573009,\Phi_2= 5.71361\times 10^{7}$} \label{fig:InitS:c}] {\includegraphics[width=0.8\textwidth]{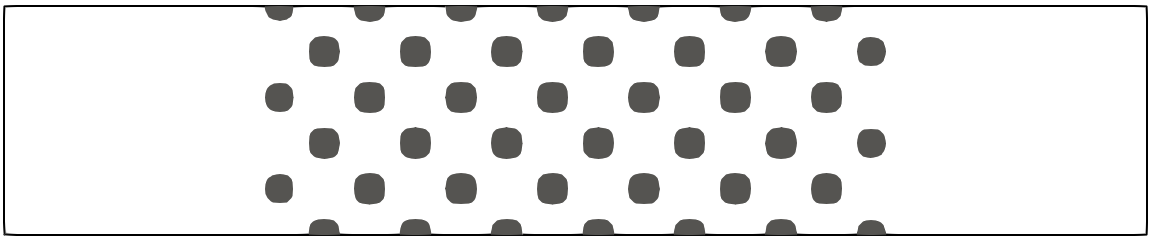}}\\
            \subfloat[{} \label{fig:InitS:d}] {\includegraphics[width=0.43\textwidth]{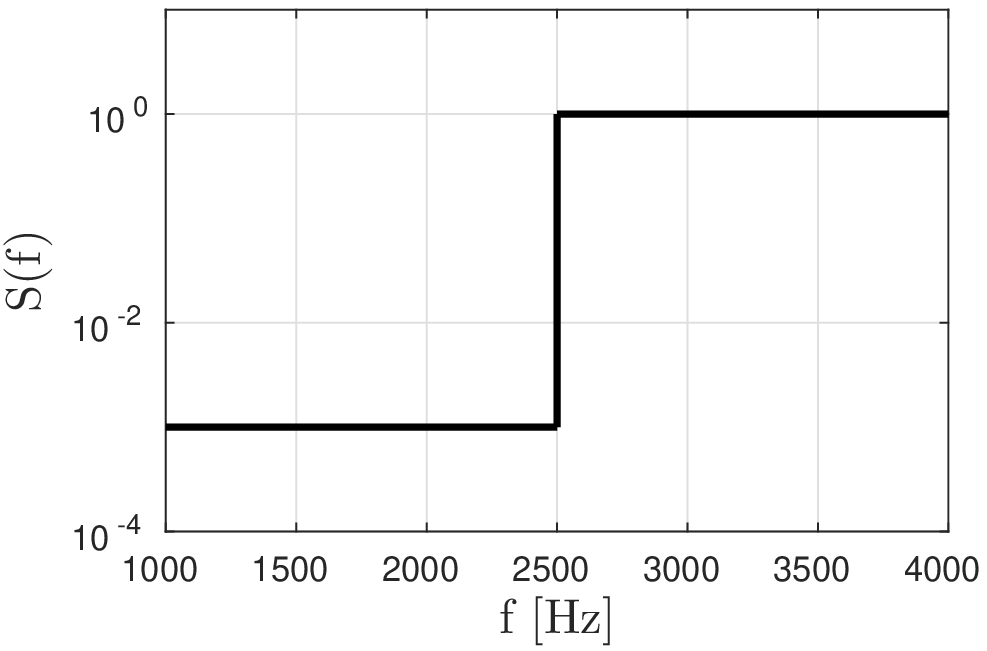}}
           \subfloat[{} \label{fig:InitS:e}] {\includegraphics[width=0.43\textwidth]{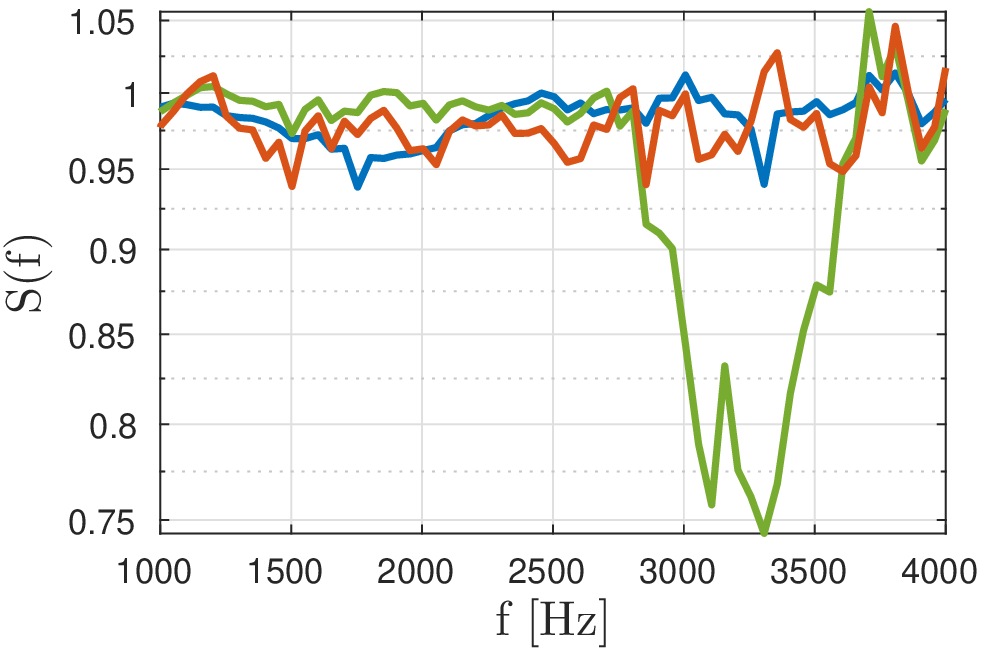}}
   \caption{Optimization for an acoustic high-pass filter design. Figures (a), (b) and (c) are the initial guess designs utilized for the study. (d) Figure showing the desired high-pass filter. (e) transmission of the initial guess designs in the considered frequency range, blue line is the response of the design shown in Fig. \ref{fig:InitS:a}, orange line is the response of the design shown in Fig. \ref{fig:InitS:b}, green line is the response of the design shown in Fig. \ref{fig:InitS:c}.}
  \label{fig:InitS}
\end{figure*}

This section carries out an initial guess study for the optimization of acoustic filters. Instead of the low-pass filter design that is presented in the previous section \ref{objStudy}, the section considers the optimization of high-pass acoustic filters. For the optimization, the \markline{constraint} function $\Phi_1$ defining the pass-band of the high-pass filter operates on the frequencies between  $2500\;\rm{Hz} \leq f \leq 4000\,\rm{Hz}$. Whereas, $\Phi_2$ for the stop-band region of the high-pass filter is active between $1000\;\rm{Hz} \leq f < 2500\,\rm{Hz}$. The desired high-pass filter is given in figure \ref{fig:InitS:d} in which the calculated transmission  $S(f)$ of the designs are to be lowered to a value $S(f) = 1 \times 10^{-3}$ in the stop-band region of the filter.

For the current study, three different initial configurations will be considered. The first initial design is the same that was used in the previous section. The other two initial guesses are obtained with decreasing the $r_1$ and $r_2$ parameters given in equation \ref{initialEquationS} from $7$ to $r_1,r_2=6$ and  $r_1,r_2=5$, respectively. Decreasing the parameters $r_1$ and $r_2$ reduce the total number of circle structures in the initial configuration while making their size bigger.

Figures \ref{fig:InitS:a} to \ref{fig:InitS:c} presents the initial configurations that are used for the optimization of acoustic high-pass filters. As it is seen from the figures, initial structures have sparsely clustered features in order to allow a high transmission across the frequencies that are considered in the high-pass filter given in figure \ref{fig:InitS:d}. Initial structures' features become larger from figure \ref{fig:InitS:a} to \ref{fig:InitS:c}. The figures also presents the \markline{constraint} values $\Phi_1$ and $\Phi_2$ calculated with the initial designs where similar results are obtained. This points out similar transmission responses between $1000\rm{Hz}$ and $4000\rm{Hz}$ for the initial configurations. 

Figure \ref{fig:InitS:e} presents the transmission $S(f)$ responses calculated at the outlet of the acoustic duct over the considered frequency range. Here, the initial configurations given in figures \ref{fig:InitS:a} and \ref{fig:InitS:b} resulted in transmission responses that closely follow each other between $1000\rm{Hz}$ and $4000\rm{Hz}$ where the calculated transmission $S(f)$ values are clustered around unity. Meaning that the initial configuration for the first two cases allows for nearly full transmission in the frequency range that is considered for the optimization. Likewise, the initial design given in figure \ref{fig:InitS:c} also allows for frequencies between $1000\rm{Hz}$ and $2800\rm{Hz}$ to pass. However, after approximately around $2800\rm{Hz}$, the last initial design has a lowered transmission response as it can also be seen from the figure \ref{fig:InitS:e} which corresponds to the pass-band region of the considered high-pass acoustic filter. 

\begin{figure*}
\centering
  \subfloat[{$\Phi_1=13.7815,\Phi_2= 11.6695$} \label{fig:InitF:a}] {\includegraphics[width=0.8\textwidth]{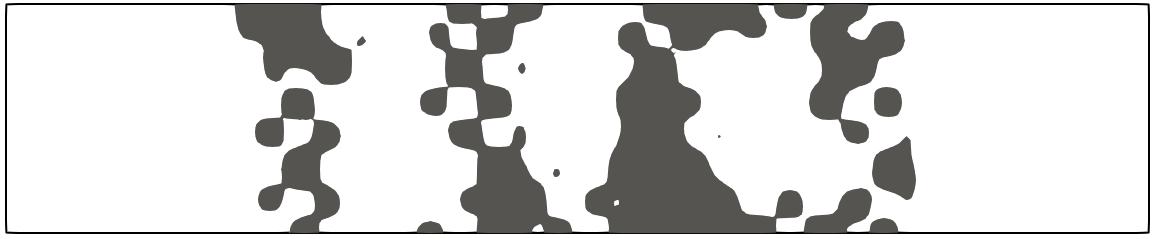}}\\
    \subfloat[{$\Phi_1=9.31977,\Phi_2= 10.0032$} \label{fig:InitF:b}] {\includegraphics[width=0.8\textwidth]{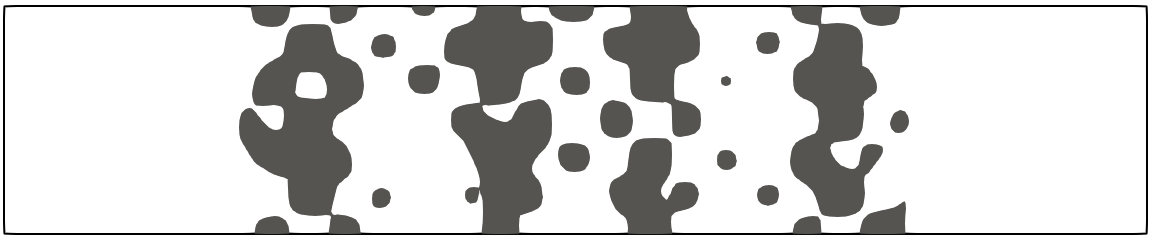}}\\
      \subfloat[{$\Phi_1=50.7457,\Phi_2= 19.9735$} \label{fig:InitF:c}] {\includegraphics[width=0.8\textwidth]{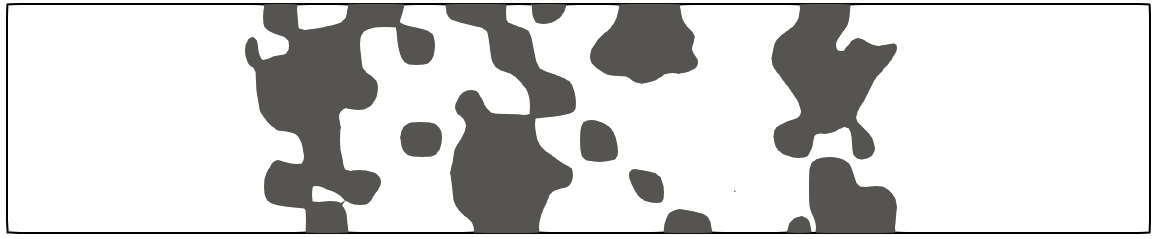}}\\
            \subfloat[{} \label{fig:InitF:d}] {\includegraphics[width=0.43\textwidth]{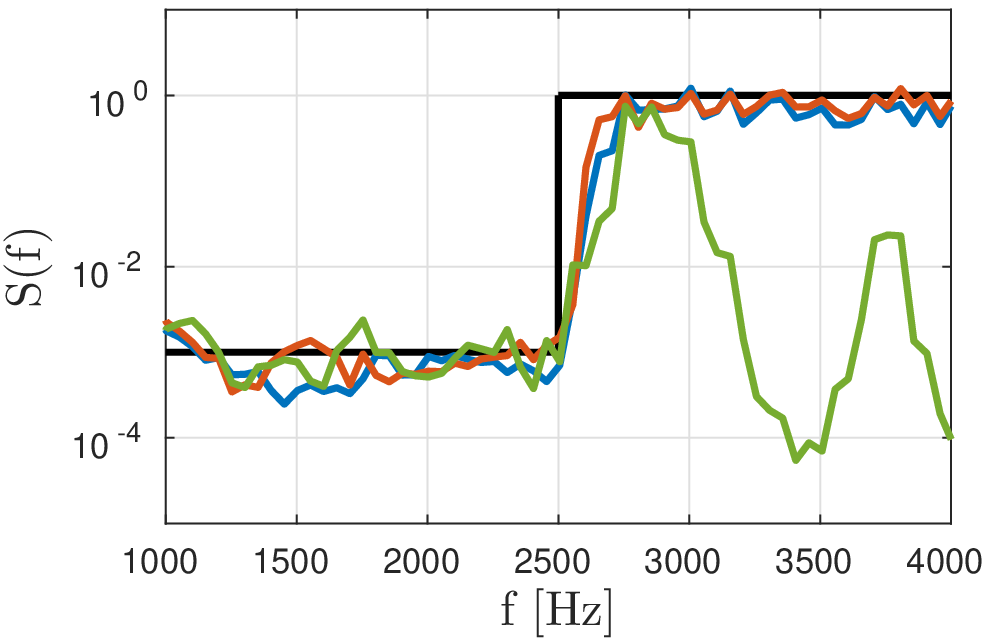}}
            \subfloat[{} \label{fig:InitF:e}] {\includegraphics[width=0.43\textwidth]{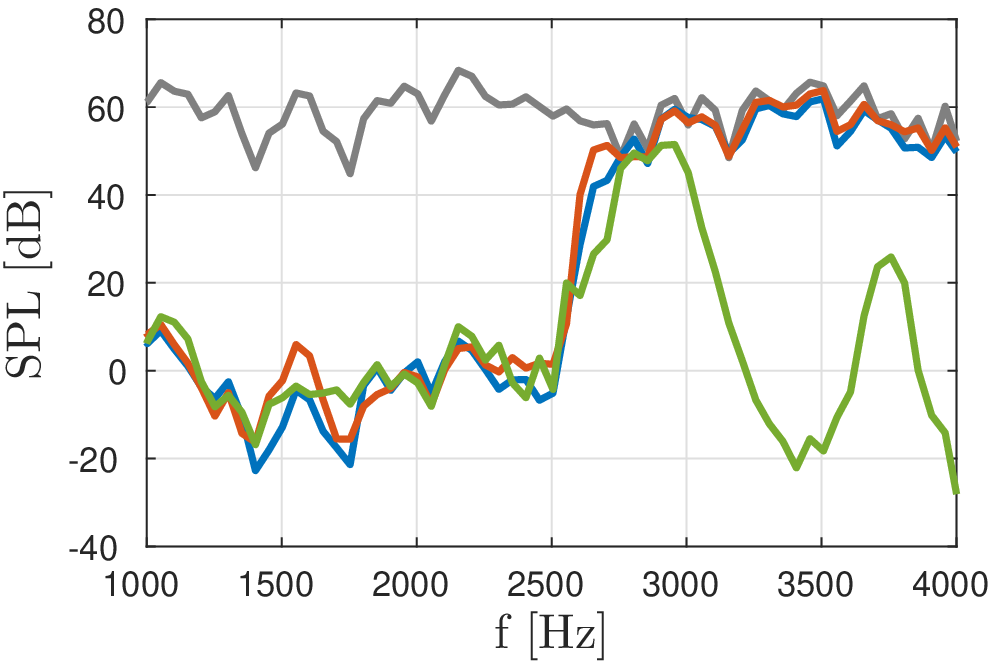}}
   \caption{Optimization for an acoustic high-pass filter design. (a) optimized design for the initial guess shown in Fig. \ref{fig:InitS:a}. (b) optimized design for the initial guess shown in Fig. \ref{fig:InitS:b}. (c) optimized design for the initial guess shown in Fig. \ref{fig:InitS:c}. (d) transmission of the optimized designs in the considered frequency range. Black line is the desired high-pass filter, blue line is the response of the design in shown Fig. \ref{fig:InitF:a}, orange line is the response of the design shown in Fig. \ref{fig:InitF:b}, green line is the response of the design shown in Fig. \ref{fig:InitF:c}. (e) averaged SPL response of the designs calculated at the output, gray line is the empty acoustic duct, blue line is the response of the design shown in Fig. \ref{fig:InitF:a}, orange line is the response of the design shown in Fig. \ref{fig:InitF:b}, green line is the response of the design shown in Fig. \ref{fig:InitF:c}. }
  \label{fig:InitF}
\end{figure*}

The results of the optimization for the initial guess study are presented in figure \ref{fig:InitF}. Moreover, the optimized designs can be seen in figures \ref{fig:InitF:a} to \ref{fig:InitF:c} which are the optimized results of the initial designs given in figures \ref{fig:InitS:a} to \ref{fig:InitS:c}, respectively. Here, the designs given in figures \ref{fig:InitF:a} and \ref{fig:InitF:b} have a similar structure layout where the optimized structures in the design domain separated the domain roughly into three acoustic partitions. \markline{Again, it is interesting to note that the optimized high pass filters do not resemble classical designs with short side branches to the main channel.} Optimized design given in figure \ref{fig:InitF:c} resulted in the largest features compared to the first two cases. Also, from the investigation of the \markline{constraint} values $\Phi_1$ and $\Phi_2$ given in figures \ref{fig:InitF:a} to \ref{fig:InitF:c}, it is seen that the first two optimized designs successfully captured the high-pass filter response whereas the last optimization failed in the pass-band region of the considered high-pass filter shape. This is also visually shown in figure \ref{fig:InitF:d} where the calculated transmission $S(f)$ for each optimized design is plotted over the considered frequency range. It can be seen in the figure that the optimized designs given in figures \ref{fig:InitF:a} and \ref{fig:InitF:b} perform successfully as acoustic high-pass filters with lowering the transmission to $S(f)= 1\times 10^{-3}$ in the stop-band and, after a sharp transition around $2500 \rm{Hz}$, realizing nearly full transmission in the pass-band region of the high-pass filter.

Furthermore, the figure \ref{fig:InitF:d} also shows the response of the design given in figure \ref{fig:InitF:c} in which the design successfully realizes the stop-band of the high-pass filter. However, after the transition into the pass-band, the design's performance deteriorates as the calculated transmission $S(f)$ lowers in the pass-band. Figure \ref{fig:InitF:e} compares the designs' averaged SPL values at the outlet of the acoustic duct for the considered frequency range against the response of an empty acoustic channel. Here for the optimized designs in figures \ref{fig:InitF:a} and \ref{fig:InitF:b}, the SPL values in the stop-band are clustered around $0 \rm{dB}$ and goes up to around $60\rm{dB}$ in the pass-band with closely following the response of an empty acoustic channel. The failure of the design given in figure \ref{fig:InitF:c} can also be seen from the figure where the SPL values diverge from the response of an empty acoustic channel in the pass-band of the high-pass filter.

Overall, the developed transient optimization framework for coupled vibroacoustic problems is successfully applied for the design of acoustic high-pass filters. Initial configuration given in figure \ref{fig:InitS:c} resulted in a failed design, \markline{which is most likely due to the fact that this initial guess  \marklineNew{does not provide full transmission for all considered frequencies. In other words, as the objective is measured on the output, the sensitivity will be zero if no signal is received at this port.} This observation is further strengthen as the} initial configurations that allow nearly full transmission (Figs. \ref{fig:InitS:a} and \ref{fig:InitS:b}) across all frequencies considered by the optimization, resulted in efficient acoustic high-pass filters.  \marklineNew{The same behavior is also seen in other works on filter design, see e.g. \cite{Aage2017ayeni} on the design of microwave waveguide filters.}

\subsection{Band-pass and band-stop acoustic filter designs}

Having validated the proposed \marklineNew{constraint} function formulation, this section concerns the design of specific filter characteristics.
The section firstly carries out the optimization for the band-pass filter where the \markline{constraint} function $\Phi_1$, defining the pass-band region of the band-pass filter, operates on the frequencies between $2500\;\rm{Hz} \leq f \leq 4000\,\rm{Hz}$. Stop-band regions defined by the \markline{constraint} function $\Phi_2$ are considered in the frequency window of  $1000\;\rm{Hz} \leq f < 2500\,\rm{Hz}$ and $4000\;\rm{Hz} < f \leq 5500\,\rm{Hz}$. Initial structure considered for the optimization is the same as in the figure \ref{fig:InitS:a} where the initial \markline{constraint} values are listed here as $\Phi_1=0.0154458,\Phi_2= 1.15859\times 10^{8}$. \markline{Remark that the two examples presented in this section are run for 800 iterations.}

\begin{figure*}
\centering
  \subfloat[{$\Phi_1=20.4625,\Phi_2= 24.4682$} \label{fig:pass:a}] {\includegraphics[width=0.8\textwidth]{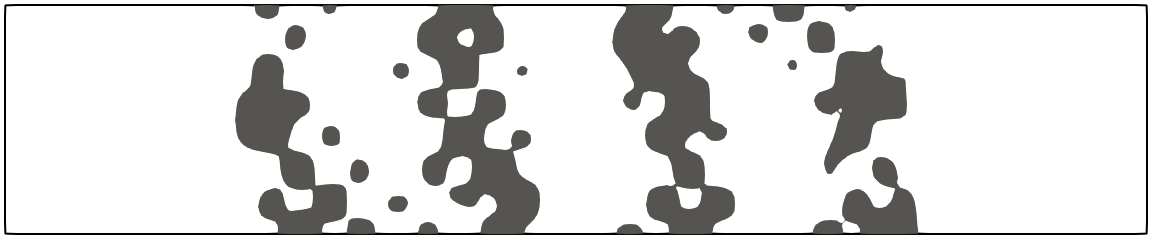}}\\
    \subfloat[{} \label{fig:pass:b}] {\includegraphics[width=0.45\textwidth]{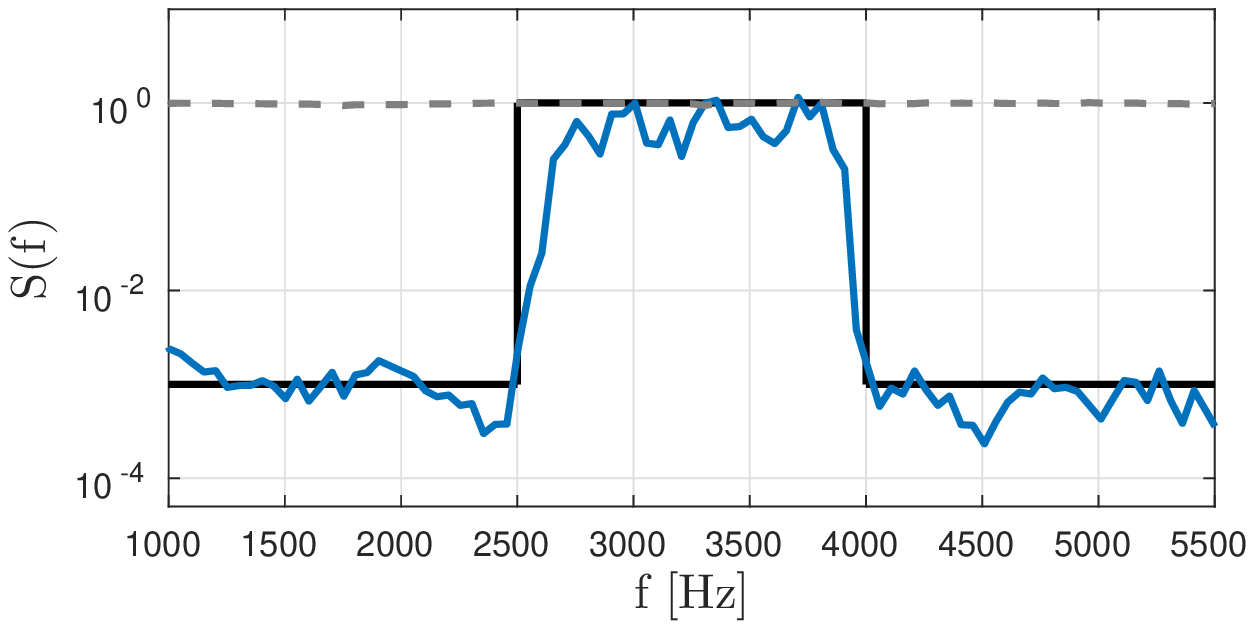}}\\
      \subfloat[{} \label{fig:pass:c}] {\includegraphics[width=0.45\textwidth]{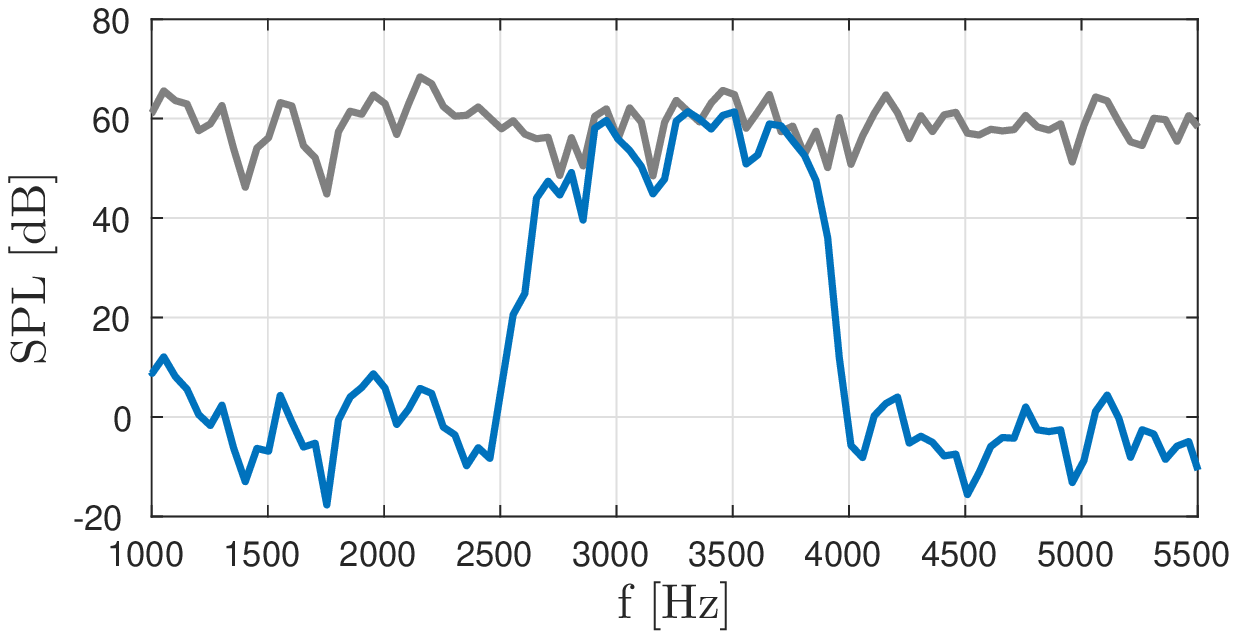}}\\
   \caption{Optimization for an acoustic band-pass filter design. (a) optimized design. (b) transmission of the optimized design in the considered frequency range. Black line is the desired band-pass filter, blue line is the response of the optimized design, gray dashed line is the response of the initial guess design. (c) averaged SPL response of the optimized design calculated at the output, gray line is the empty acoustic duct, blue line is the response of the optimized design.}
  \label{fig:pass}
\end{figure*}

Figure \ref{fig:pass:a} shows the optimized \markline{band-pass filter} while the transmissions $S(f)$ of the optimized design and the initial configuration are given in figure \ref{fig:pass:b}. As it can be seen from  figure \ref{fig:pass:b}, the initial configuration has nearly full transmission across the frequencies $1000\;\rm{Hz}$ to $5500\,\rm{Hz}$. The optimized design on the other hand closely follows the desired band-pass filter shape. Interestingly, the optimized band-pass filter (Fig.  \ref{fig:pass:a}) has a similar structural arrangement compared to the high-pass filter in figure \ref{fig:InitF:a} and in which both designs have similar transmission responses from $1000\;\rm{Hz}$ to $4000\,\rm{Hz}$ for a high-pass filter shape. 
 The current design further tailors the frequency response of the coupled system to realize an additional stop-band from $4000\;\rm{Hz}$ to $5500\,\rm{Hz}$. Also, from the optimized design's \markline{constraint} values given in figure \ref{fig:pass:a} and the visual inspection of the transmission response in figure \ref{fig:pass:b}, it can be said that the performance of the pass-band region is slightly lowered compared to the designs presented in previous sections. This is mainly due to the added complexity of the optimization where an additional stop-band is considered to realize a band-pass filter. 

Figure \ref{fig:pass:c} shows \marklineNew{the calculated averaged SPL response at the outlet} of the acoustic channel, comparing the optimized design and the empty acoustic channel. It can be seen from the figure that the optimized design's calculated SPL response at the outlet is clustered around $0\rm{dB}$ for the stop-band regions of the band-pass filter. At the pass-band, the design's SPL response closely follow the SPL response of an empty acoustic channel. Overall, the developed transient optimization framework is shown to successfully tailor the frequency response of the coupled vibroacoustic system with designing an efficient acoustic band-pass filter.

The framework is lastly applied for the design of a band-stop filter. The \markline{constraint} function $\Phi_1$ for the pass-band regions of the considered band-stop filter operates on the frequency window of $1000\;\rm{Hz} \leq f < 2500\,\rm{Hz}$ and $4000\;\rm{Hz} < f \leq 5500\,\rm{Hz}$. For the stop-band region, the \markline{constraint} function $\Phi_2$ is active between the frequencies of $2500\;\rm{Hz} \leq f \leq 4000\,\rm{Hz}$. Again, the same initial configuration as the previous band-pass filter optimization is utilized which resulted in initial \markline{constraint} values of $\Phi_1=0.0552177,\Phi_2= 5.8754\times 10^{7}$.

\begin{figure*}
\centering
  \subfloat[{$\Phi_1=24.8701,\Phi_2= 60.3181$} \label{fig:stop:a}] {\includegraphics[width=0.8\textwidth]{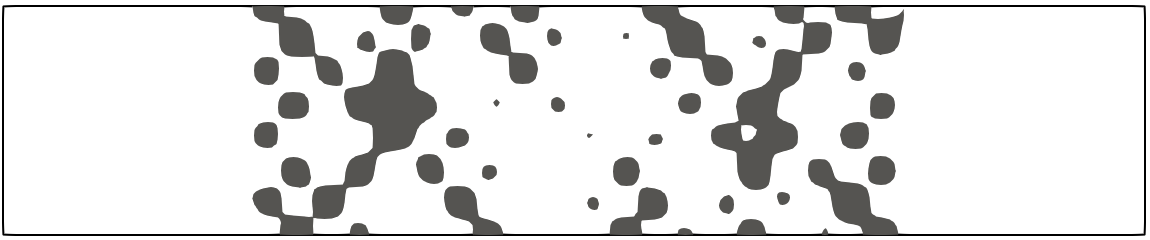}}\\
    \subfloat[{} \label{fig:stop:b}] {\includegraphics[width=0.45\textwidth]{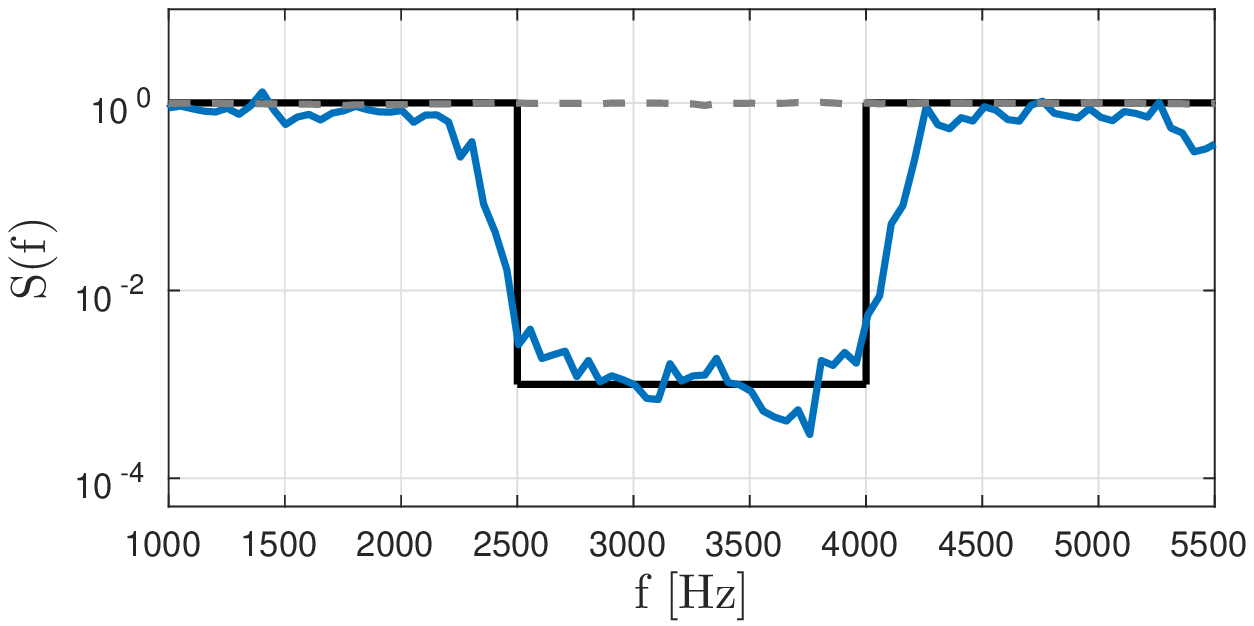}}\\
      \subfloat[{} \label{fig:stop:c}] {\includegraphics[width=0.45\textwidth]{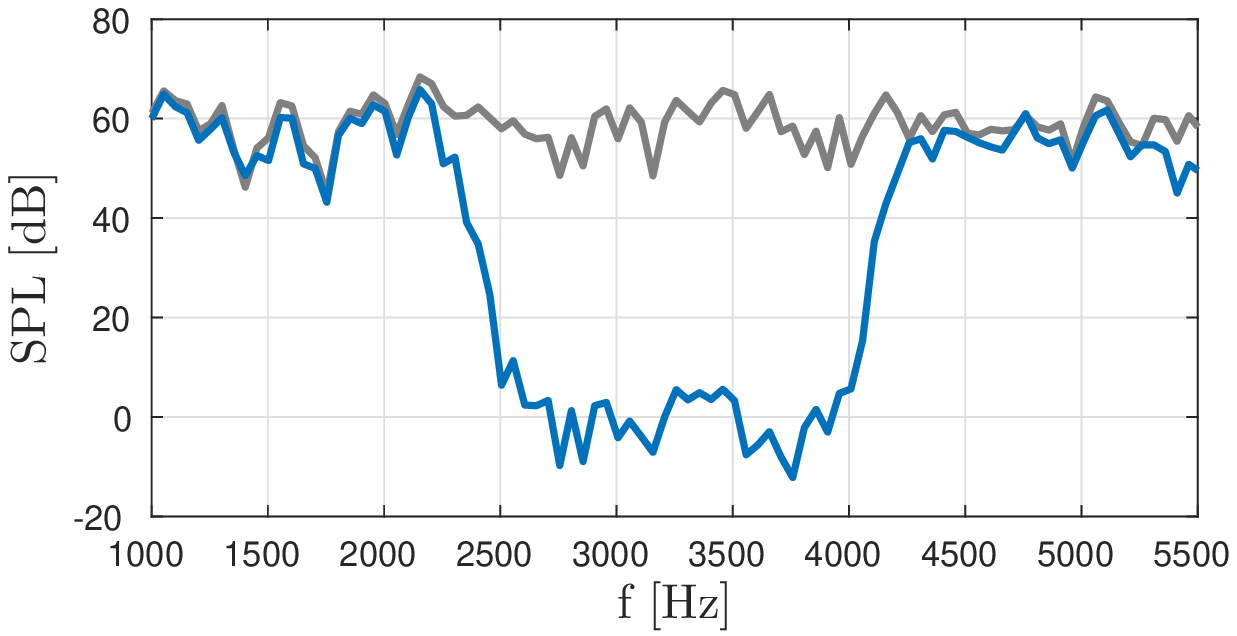}}\\
   \caption{Optimization for an acoustic band-stop filter design. (a) optimized design. (b) transmission of the optimized design in the considered frequency range. Black line is the desired band-pass filter, blue line is the response of the optimized design, gray dashed line is the response of the initial guess design. (c) averaged SPL response of the optimized design calculated at the output, gray line is the empty acoustic duct, blue line is the response of the optimized design.}
  \label{fig:stop}
\end{figure*}

The optimized design is shown in figure \ref{fig:stop:a}. Moreover, the desired band-stop filter shape along with the transmissions $S(f)$ of the optimized and initial designs are given in figure \ref{fig:stop:b}. When the optimized design is compared to the low-pass filter design given in figure \ref{fig:AB:b}, it can be said that the both designs have a similar performance between the frequencies from $1000\;\rm{Hz}$ to $4000\,\rm{Hz}$ with acting as an acoustic low-pass filter. Overall, it can be seen in figure \ref{fig:stop:b} that the initial transmission $S(f)$ is successfully tailored during the optimization in which the optimized design's response closely follow the desired band-stop filter between the frequencies of $1000\;\rm{Hz}$ to $5500\,\rm{Hz}$. The performance of the optimized acoustic band-stop filter is also seen from the averaged SPL values calculated at the outlet of the acoustic channel which is given in figure  \ref{fig:stop:c}. The figure shows that pass-band regions attains the overall $60\rm{dB}$ SPL response while the stop-band region of the filter lowers the SPL response to approximately around $0\rm{dB}$.

\markline{As was the case with the previous optimization results, the optimized band-pass and band-stop filters share little similarity to classical filter design, i.e. no clearly distinguishable Helmholtz resonators, narrow slits, channel extensions are seen. }

\subsection{Optimization and validation of an acoustic low-pass filter} 
\label{sec:validation}
This section concerns the optimization of an acoustic low-pass filter and aims
to verify the performance of the optimized design using a commercial software. 
In order to achieve this, after the optimization is finished, the obtained design is extracted and \markline{both solid and acoustic regions are re-}meshed using body-fitted\markline{, conforming meshes consisting of a total of 151460 triangles with quadratic shape functions}. \marklineNew{This means that even the smallest islands seen in the optimized designs are meshed with several elements in order to ensure proper capturing of the physics. It should also be noted that, the optimized design is directly obtained from the the physical design variable $\bar{s}$ and no smoothing operation has been applied on the interface before re-meshing.} The widely used commercial software COMSOL Multiphysics \cite{comsolRef} is used \markline{for the validation}, in which the strongly coupled vibroacoustic problem is solved \markline{on time-harmonic form.
Thus, the displacement and acoustic pressure degrees of freedoms in the body-fitted mesh are only defined in their corresponding solid and acoustic domains and only share the coupled interface which is explicitly defined. That is, COMSOL employs no ersatz material model or special integration rules as is the case with the cut element method. Moreover, we should emphasize strongly that using time-harmonics, i.e. a steady-state formulation, to post evaluate the performance of the optimized designs which are obtained using a transient formulation, is a very hard test to pass. That is, the transient optimization problem implementation is likely to utilize any numerical and/or physical artifact in the underlying model in order to gain performance improvements. For the problems at hand, this includes differences in loading, exploitation of transients, constructive/destructive interference, dissipation in the time-stepping scheme, etc. None of these effects will be present in the time-harmonic formulation. }  

As it has been thoroughly introduced in previous sections, the transmission $S(f)$ of the optimized design is calculated with applying the FFT operation on the transient response of the coupled system when the system is excited with the incoming white noise realized as pressure oscillations. In order to compare and validate the calculated transmission $S(f)$ of the design against the body-fitted mesh analysis done in COMSOL, we utilize a time-harmonic frequency domain analysis to achieve the frequency-sweep where the system is excited with an incoming sinusoidal plane wave for each frequency. This means that the system is solved at discrete frequencies within the frequency range that is considered by the optimization to realize steady-state solutions at each frequency. The transmission $S(f)$ is then calculated at each discrete frequency and compared against the one obtained with the developed transient framework. \markline{The rest of the boundary conditions considered in COMSOL are; both top and bottom boundaries of the acoustic channel are set as hard-wall condition for the acoustic pressure while the structure is considered to be clamped, the outlet of the acoustic channel is set as an absorbing boundary and the coupling boundary conditions are applied to both acoustic pressure and displacement variables at the coupled interface.}

The acoustic low-pass filter is considered for the frequency range of $500\;\rm{Hz} \leq f \leq 3500\,\rm{Hz}$ where the pass-band acts on the frequencies $500\;\rm{Hz} \leq f \leq 2000\,\rm{Hz}$ and the stop-band is defined on the window $2000\;\rm{Hz} < f \leq 3500\,\rm{Hz}$. The desired transmission for the stop-band is considered as $S(f)= 1\times 10^{-4}$ to be able to realize in average four orders of magnitude decrease in the amplitude of the transmitted acoustic pressure compared to that of an empty channel.

\begin{figure*}
\centering
  \subfloat[{$\Phi_1=14.4407,\Phi_2= 79.0555$} \label{fig:compa:a}] {\includegraphics[width=0.8\textwidth]{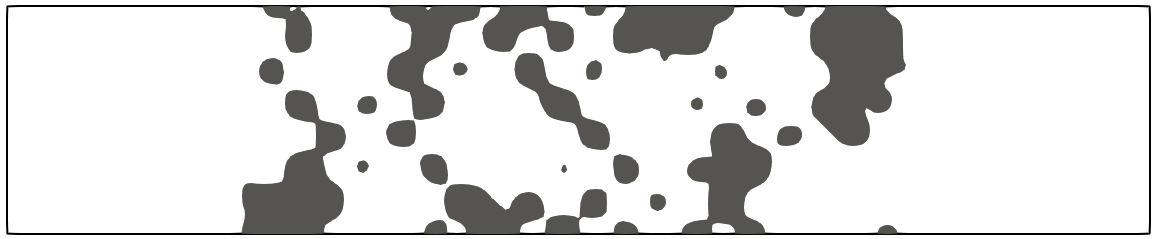}}\\
      \subfloat[{} \label{fig:compa:b}] {\includegraphics[width=0.45\textwidth]{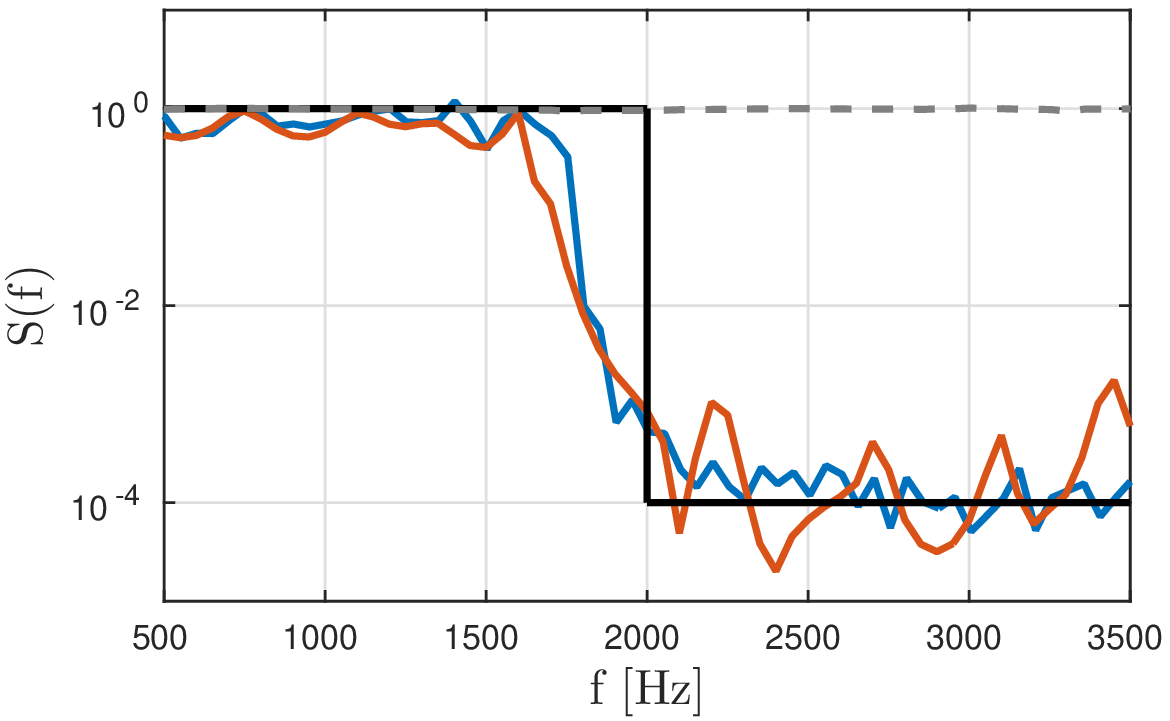}}\\
   \caption{Optimization for an acoustic low-pass filter design. (a) optimized design. (b) transmission of the optimized design in the considered frequency range. Black line is the desired low-pass filter, blue line is the response of the optimized design, gray dashed line is the response of the initial guess design and the orange line is the COMSOL calculation done in frequency domain with body-fitted analysis.}
  \label{fig:compa}
\end{figure*}

For optimization, the initial structure given in figure \ref{fig:ABinit:a} is used which results in initial \markline{constraint} values of $\Phi_1=0.0385819,\Phi_2= 5.8394\times 10^9$. These values reflect nearly full transmission in the considered frequency range which can also be seen from figure \ref{fig:compa:b} where the response of the initial guess design is plotted on the desired filter shape.

The result of the acoustic low-pass filter optimization is presented in figure \ref{fig:compa} and  figure \ref{fig:compa:a} shows the optimized design. Similarly to the previous optimization studies,  it can be seen that the optimized design resulted in significant topological changes compared to the initial configuration\markline{, i.e. several solid regions has merged and some has completely} \marklineNew{disappeared.} Moreover, figure \ref{fig:compa:b} presents the calculated transmission $S(f)$ of the optimized design along with the transmission of the post-processed design obtained with frequency domain calculation using COMSOL. As it can be seen from the figure, the calculated transmission \markline{response obtained} using the transient framework with the cut element method and the body-fitted analysis done in frequency domain with COMSOL agrees well. For both \markline{analysis} methods, the calculated transmission \markline{response} closely follows the desired pass-band and stop-band regions. As expected, the transition from the pass-band into the stop-band (around $2000\,\rm{Hz}$) causes a slight decrease in the performance of the pass-band region. However, the design successfully lowers the transmission to approximately around $S(f)= 1\times 10^{-4}$ at the stop-band.
\markline{Despite the fact that the two analysis method cannot be directly compared, we here quantify the observed mismatch using an average error measure. That is, the absolute difference in the average transmission $<S(f)>$ between the transient cut element analysis and the body-fitted mesh time-harmonic COMSOL analysis is computed in the pass-band and stop-band and denoted with $(\cdot)_\text{cut}$ and $(\cdot)_\text{body}$, respectively.
\begin{align*}
    &|<S(f)>_\text{cut} - <S(f)>_\text{body}|_{\text{pass}}=0.107\\
    &|<S(f)>_\text{cut} - <S(f)>_\text{body}|_{\text{stop}}=2.237\cdot 10^{-4}
\end{align*}}
\markline{where one should keep in mind that the pass-band has unity as perfect transmission whereas the stop-band should have a transmission of $S(f)= 1\times 10^{-4}$. Returning to the inspection of figure \ref{fig:compa:b}, we remark that even though the two analysis are performed with quite different modelling and analysis strategies, both methods leads to qualitatively similar transmission responses. To elaborate,} the discrepancies in the responses of the transient and frequency domain methods can be explained from the fact that, the frequency domain analysis results in steady-state solutions at each discrete frequency whereas the developed framework utilizes a finite transient signal to calculate the transmission through the FFT operation. Since the transient system is excited with random pressure oscillations it would be computationally unreasonable to consider steady-state solutions as this would be too time consuming. \markline{The difference in modelling and analysis is therefore attributed the majority of the difference in response.} Another factor that can cause discrepancies in the response when comparing two method is the meshing procedure in which when the design is extracted from the zero iso-level of $\bar{s}$ and meshed with unstructured meshes, a slight alteration of the interface is unavoidable. Overall, it can be said that, when analyzed with more commonly used methods such as time-harmonic frequency domain methods using body-fitted analysis, the optimized design also shows a good performance as an acoustic low-pass filter.

\begin{figure*}
\centering
 \subfloat{\includegraphics[width=0.5\textwidth]{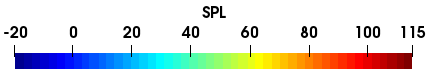}} \addtocounter{subfigure}{-1}  \\
  \subfloat[{} \label{fig:compa2:a}] {\includegraphics[width=0.8\textwidth]{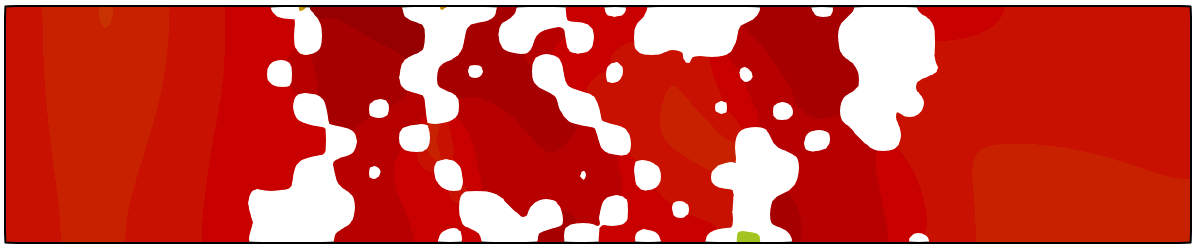}}\\
  \subfloat[{} \label{fig:compa2:b}] {\includegraphics[width=0.8\textwidth]{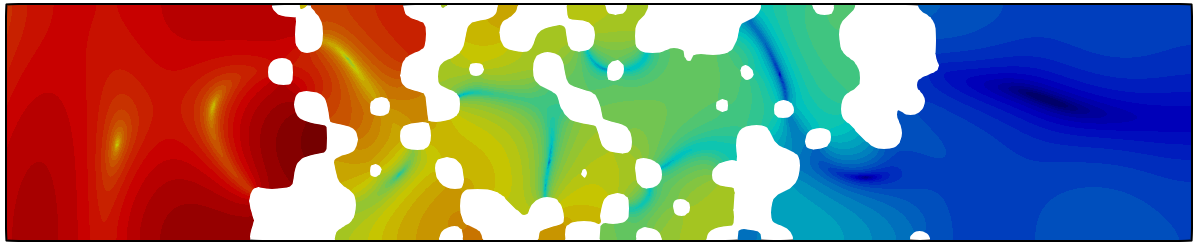}}\\
   \caption{\markline{COMSOL calculation of} time-harmonic frequency domain analysis of the optimized design for the acoustic low-pass filter showing the sound pressure level $\text{[dB]}$ contours of the design. (a) The \markline{steady state} result from the pass-band at $f=1100\,\text{Hz}$. (b) The \markline{steady state} result from the stop-band at $f=2900\,\text{Hz}$. }
  \label{fig:compa2}
\end{figure*}

\begin{figure*}
\centering
 \subfloat{\includegraphics[width=0.5\textwidth]{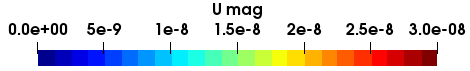}} \addtocounter{subfigure}{-1}  \\
  \subfloat[{} \label{fig:compa3:a}] {\includegraphics[width=0.8\textwidth]{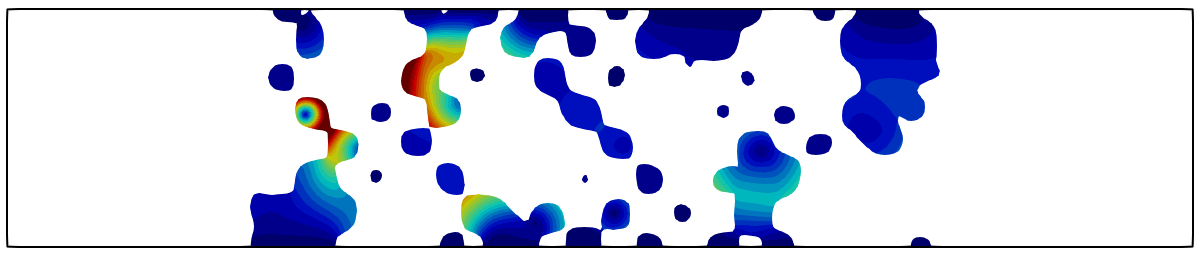}}\\
  \subfloat[{} \label{fig:compa3:b}] {\includegraphics[width=0.8\textwidth]{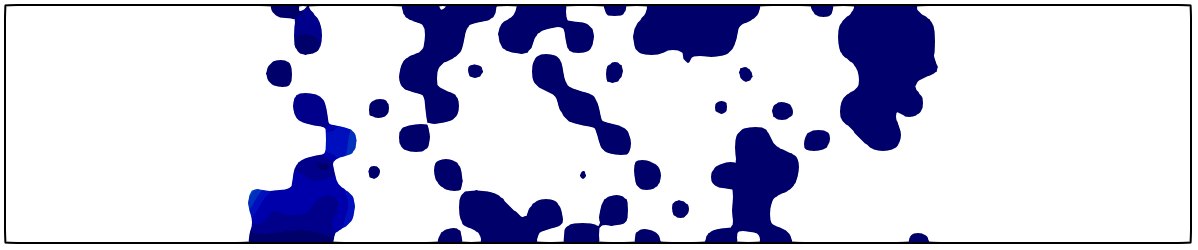}}\\
   \caption{\markline{COMSOL calculation of} time-harmonic frequency domain analysis of the optimized design for the acoustic low-pass filter showing the displacement magnitude $|\mathbf{u}|$ $\text{[m]}$ contours of the design. (a) The \markline{steady state} result from the pass-band at $f=1100\,\text{Hz}$. (b) The \markline{steady state} result from the stop-band at $f=2900\,\text{Hz}$. }
  \label{fig:compa3}
\end{figure*}

Figures \ref{fig:compa2} and \ref{fig:compa3} present the corresponding SPL and displacement magnitude, respectively, obtained from the time-harmonic analysis of the optimized design using COMSOL. The figures consider two discrete frequencies $f=1100\,\text{Hz}$ and $f=2900\,\text{Hz}$ for pass-band and stop-band, respectively. Furthermore, figure \ref{fig:compa2:a} shows the SPL field of the optimized design where the pressure response is plotted at $f=1100\,\text{Hz}$ (at the pass-band). As it can be seen from the figure, the design allows the incoming wave to pass resulting in approximately constant SPL response at around $110\rm{dB}$. Corresponding displacement magnitude $|\mathbf{u}|$ $\text{[m]}$ field is plotted in figure \ref{fig:compa3:a} in which complex deflections of the structure can be seen which means that the underlying acoustic-structural interactions are fully utilized to realize the desired transmission behaviour of the design. On the other hand when the design is analysed at $f=2900\,\text{Hz}$ (at the stop-band), as seen from figure \ref{fig:compa2:b}, the incoming wave is stopped and the design realizes SPL values of below 0 dB at the outlet of the acoustic channel. Interestingly, as it is seen from the corresponding displacement magnitude $|\mathbf{u}|$ $\text{[m]}$ field in figure \ref{fig:compa3:b}, the resulted deflections of the structure are at least one order of magnitude smaller than the solution presented at the pass-band. Meaning that the optimized design acts as a wave stopper primarily because of the internal reflections and the underlying acoustic-structural interactions play a smaller part.

\begin{figure*}
\centering
  \subfloat[{} \label{fig:desMod:a}] {\includegraphics[width=0.9\textwidth]{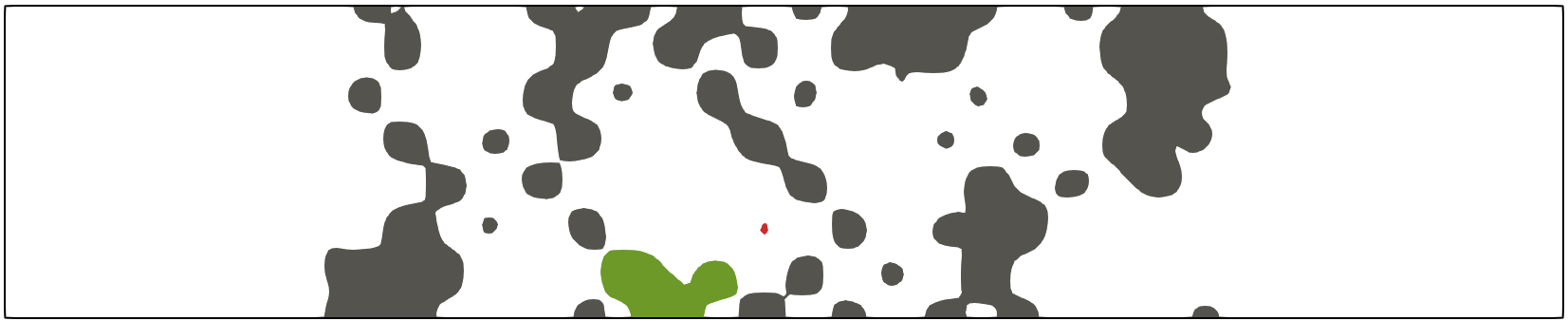}}\\
      \subfloat[{} \label{fig:desMod:b}] {\includegraphics[width=0.5\textwidth]{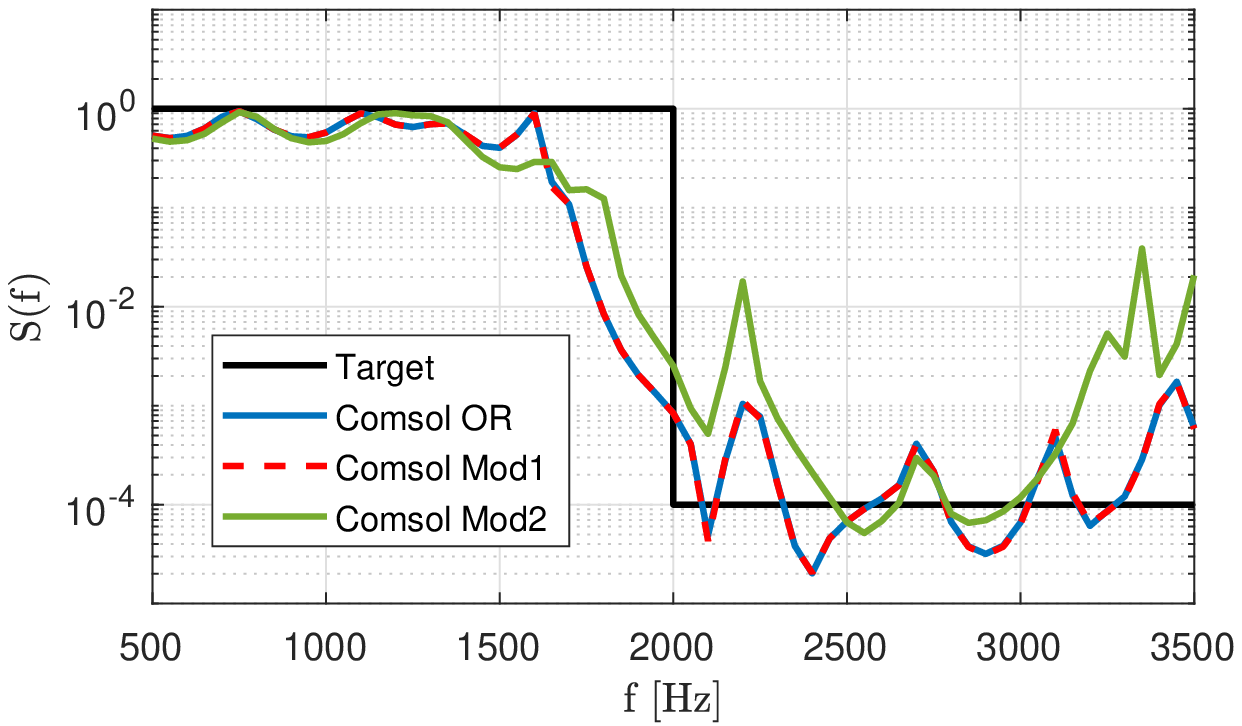}}\\
   \caption{\marklineNewL{(a) Figure shows the locations of the modifications to the design. The first modification is the red colored island is removed and the rest of the structure (including the green colored part) is analyzed. The second one is the green colored part is removed and the rest of the structure (including the red colored part) is analyzed. (b) Transmission of the structure in the considered frequency range where black line is the target low-pass filter, blue line is the response of the original design, red dashed line is the response of the first modified design and the green line is the response of the second modification. All calculations have been done using COMSOL in frequency domain with body-fitted analysis.}}
  \label{fig:desMod}
\end{figure*}

\marklineNewL{Lastly, the optimized structure is further analysed to asses the effects of small island inclusions on the overall performance of the design in terms of its low-pass filter capabilities. In order to carry out the study, two island parts of the design are selected where the first one is the smallest and the other is relatively bigger part of the overall design. Figure \ref{fig:desMod:a} shows the selected regions where the each part is taken out in turn and the design is analysed in time-harmonic frequency domain analysis using COMSOL for the same selected frequency range of $500\,\text{Hz}$ to $3500\,\text{Hz}$ while comparing the results to the COMSOL
calculation of the original design. The result of the comparative study is given in figure \ref{fig:desMod:b} in which it can be seen that excluding the smallest island part of the design does not change the performance of the design. This can be explained from the fact that at the frequency range of $500\,\text{Hz}$ to $3500\,\text{Hz}$, the wavelength of air is significantly larger than the size of the considered part. Thus, it does not make any difference if the part is included in the design or not. Unsurprisingly, excluding the second part from the design results in significant performance reduction especially in the stop-band region since the internal reflections of the channel gets disturbed which results in lowered performance.}

\marklineNewL{These extremely small island structures (like the red colored part shown in figure \ref{fig:desMod:a}) also have been consistently present in the previous optimization cases. The lack of volume constraint together with the highly non-convex nature of the optimization problem can be considered as the biggest factors of these features since the optimization becomes unbounded in the amount of material that can be used while converging to a local minima. They can also be attributed to the fact that the utilized transient formulation with the FFT operation is prone to numerical noise which may generate numerical and/or physical artifacts throughout the optimization. Once these small parts occur, the sensitivities are likely extremely small on their interfaces which can make it hard for the optimizer to remove or move them within the optimization iterations. It is important to note that these small features can also adversely affect the oscillatory behaviour of the objective function history (see figure \ref{fig:CONGRAPH}). Using a tight volume constraint with a larger smoothing filter on the design variable can help to alleviate this problem. Utilization of minimum length scales using e.g. geometric constraints \cite{Zhou2015} and/or the robust design approach \cite{wang2011a,Schousboe2020} should also be considered for realizing designs for manufacturing.}

\section{Discussion and conclusion}

The article presents the utilization of time-domain methods to realize wideband optimization in frequency domain. The developed framework carries out generalized shape optimization of transient vibroacoustic problems in which the optimization considers acoustic filter designs. In order to achieve this, the objective and constraint functions are defined in frequency domain where the FFT operation is utilized to obtain the frequency response from the transient response of the coupled system. Throughout the work, the level set approach is utilized for the geometry description where its zero iso-level specifies the interface between acoustic and structural domains. An immersed boundary method, i.e. the cut element method, is employed for capturing the geometry which operates on a fixed background mesh. The method uses a special integration scheme to accurately resolve the interface between the two physics without the addition of extra degrees of freedom to the system. Hence, the employed cut element method is suitable to include into the existing parallel FEM frameworks with ease. Moreover, the work utilizes the discrete adjoint method for carrying out the sensitivity analysis in order to calculate the gradients
of the constraint functions. The derivation of the sensitivity analysis is kept general to allow for the inclusion of different time integration schemes in which the handling of the FFT operation is explained to be able to define objective and constraint functions in frequency domain. Furthermore, a study for the utilized \markline{optimization formulation} \marklineNew{is carried out to assess the effect} of the inverse weighting that is used in the \markline{constraint} function where a design of acoustic low-pass filter is considered. It has been found that having $b=1\times 10^{-3}$ for the \markline{constraint} function $\Phi_2$ to realize the stop-band region of the considered filter resulted in a relatively sharp transition between the pass-band and the stop-band, realizing an effective filter performance. The selected $b$ parameter also reduces the averaged SPL values at the outlet of the acoustic channel to approximately around $0\rm{dB}$ in the stop-band. An initial guess study is then carried out the determine the effect of different initial configurations to the end design considering the optimization of an acoustic high-pass filters. It has been illustrated that with using initial configurations that allow for nearly full transmission in the frequency range that is considered for the optimization, efficient acoustic high-pass designs are obtained. The outcomes from both the \markline{constraint} functions and the initial guess studies are then utilized in optimization for more complex acoustic band-pass and band-stop filters. Here it is shown that the proposed optimization formulation is capable of producing filters with up to 60 dB difference between pass-bands and stop-bands. Furthermore, the result of the \markline{time-harmonic, steady-state COMSOL} validation study \markline{demonstrate} the applicability of the developed transient framework for broadband applications considering coupled vibroacoustic systems. Overall, the developed transient optimization framework is shown to successfully tailor the frequency response of the coupled vibroacoustic system for the design of acoustic band-pass and band-stop filters. The presented work will pave the way for the optimization of acoustic devices to realize efficient wideband operation. 

\markline{The proposed optimization framework presents several interesting directions which should be investigated in greater detail. The choice of constraint functions, especially their scaling, should be further examined. Although the proposed constraint functions leads to working filters with controllable characteristics, it would be desirable to determine constraint functions that results in less oscillations over the duration of the optimization process. One could also consider continuation schemes to remedy this issue, e.g. schemes in which the design variable move limits are tightened as the optimization progresses. This would also allow to include a better stopping criteria based on the constraint function evolution and the physical design change. It would also be interesting to expand the design parameterization to include manufacturing requirements. This includes imposition of minimum length scales using e.g. geometric constraints \cite{Zhou2015} or the robust design approach \cite{wang2011a,Schousboe2020}. 
Moreover, extension to 3D and treatment of the freely flying structures must also be addressed before manufacturing and application to industrially relevant problems. The latter could be done through the inclusion of additional connectivity constraints such as non-zero structural eigenfrequency requirements, the virtual temperature method or similar approaches.
}

\section*{Replication of results}
For replicating the presented examples, all necessary information is given in the corresponding sections. Moreover, the code can be obtained from the authors upon reasonable request.

\section*{Compliance with ethical standards}

\begin{description}[leftmargin=0cm]
\item[\small\textbf{Conflict of interest}]\small The authors declare that they have no conflict of interest.\small \bigskip
\item[\small\textbf{Funding}]\small There is no funding source.\small \bigskip
\item[\small\textbf{Ethical approval}]\small This article does not contain any studies with human participants or animals performed by any of the authors.
\end{description}

\bibliographystyle{spmpsci}      
\bibliography{mybibfile}   

\end{document}